\documentclass[11pt,a4paper]{article}
\usepackage{ifthen}
\usepackage{amsmath}
\usepackage{amssymb}
\usepackage{latexsym}
\usepackage{bbm}
\usepackage{graphicx}
\usepackage{pdfsync}
\usepackage{verbatim}
\usepackage{fancyhdr}
\usepackage{url}               
\usepackage[english]{babel}
\usepackage[applemac]{inputenc}
\usepackage{ mathrsfs }

\usepackage[pdftex,%
a4paper=true,%
backref=true,%
pagebackref=true,%
bookmarks=true,%
bookmarksnumbered=true,%
colorlinks=true,%
urlcolor=blue,%
pdfstartview=FitH%
]{hyperref}

\newcommand{\be}{\begin{equation}}
\newcommand{\ee}{\end{equation}}
\newcommand{\ba}{\begin{aligned}}
\newcommand{\ea}{\end{aligned}}
\newcommand{\bp}{{\it Proof. }}
\newcommand{\ep}{\hfill $\square$\\}
\newcommand{\N}{\mathbb N}
\newcommand{\NN}{\mathcal N}

\newcommand{\G}{\mathcal G}
\renewcommand{\O}{\mathcal O}
\newcommand{\R}{\mathbb R}

\newcommand{\F}{\mathcal F}

\newcommand{\A}{\mathcal A}
\newcommand{\Y}{\mathcal Y}

\newcommand{\X}{\mathcal X}
\newcommand{\tY}{\widetilde Y}

\newcommand{\tT}{\widetilde T}
\newcommand{\tv}{\vec t}
\newcommand{\bg}{ \vert\negthinspace\negthinspace  \vert}

\def\<{\langle}
\def\>{\rangle}
\renewcommand{\r}{\mathfrak{r}}
\newcommand{\supp}{\operatorname{supp}}

\newtheorem{lem}{Lemma}[section]
\newtheorem{thm}[lem]{Theorem}
\newtheorem{prop}[lem]{Proposition}
\newtheorem{cor}[lem]{Corollary}

\newtheorem{exa}[lem]{Example}

\newtheorem{remk}[lem]{Remark}

\numberwithin{equation}{section}
\numberwithin{figure}{section}

\begin{document}

\renewcommand{\refname}{References}
\bibliographystyle{alpha}

\pagestyle{fancy}
\fancyhead[L]{ }
\fancyhead[R]{}
\fancyfoot[C]{}
\fancyfoot[L]{ }
\fancyfoot[R]{}
\renewcommand{\headrulewidth}{0pt} 
\renewcommand{\footrulewidth}{0pt}

\newcommand{\montitre}{Spectral analysis of hypoelliptic random walks }

\newcommand{\auteur}{\textsc{ Gilles Lebeau, Laurent Michel}}
\newcommand{\affiliation}{Laboratoire J.-A. Dieudonn\'e \\
 Universit\'e de Nice Sophia-Antipolis\\
 Parc Valrose, 06108 Nice Cedex 02,   France\\
\url{lebeau@unice.fr, lmichel@unice.fr}}

 \begin{center}
{\bf  {\LARGE \montitre}}\\ \bigskip \bigskip
 {\large\auteur}\\ \bigskip \smallskip
 \affiliation \\ \bigskip
\today
 \end{center}

\begin{abstract}
We study the spectral theory of a reversible Markov chain 
associated to a hypoelliptic random walk on a manifold $M$. This random walk
 depends on a parameter $h\in ]0,h_{0}]$ which is roughly the size of each 
step of the walk. We prove
 uniform bounds with respect to $h$ on the rate of convergence to equilibrium,
 and the convergence when $h\rightarrow 0$ to the associated hypoelliptic diffusion.
\end{abstract}

\tableofcontents
\newpage

 \pagestyle{fancy}
\fancyhead[R]{\thepage}
\fancyfoot[C]{}
\fancyfoot[L]{}
\fancyfoot[R]{}
\renewcommand{\headrulewidth}{0.2pt} 
\renewcommand{\footrulewidth}{0pt} 

\section{Introduction and Results}\label{sec0}

The purpose of this paper is to study the spectral theory of a reversible Markov chain 
associated to a hypoelliptic random walk on a manifold $M$. This random walk
will depend on a parameter $h\in ]0,h_{0}]$ which is roughly the size of each 
step of the walk. We are in particular interested, as in \cite{DiLeMi11} and \cite{DiLeMi11},
to get uniform bounds with respect to $h$, on the rate of convergence to equilibrium.
The main tool in our approach is to compare the random walk on $M$ with a natural random walk on a nilpotent Lie group. 
This idea was used by Rotschild-Stein \cite{RoSt76} to prove sharp hypoelliptic  estimates for some differential operators.\\

Let $M$ be a smooth, connected, compact  manifold of dimension $m$, equipped with a smooth volume form
$d\mu$ such that $\int_{M} d\mu =1$. We denote by $\mu$ the associated probability on $M$. Let $\X=\{X_1,\ldots, X_p\}$ be a collection of  smooth vector fields on $M$. Denote $\G$ the Lie algebra generated by $\X$.
In all the paper we assume that the $X_k$ are divergence free with respect to $d\mu$
\be\label{hyp_nulldiv}
\forall k=1,\ldots,p, \quad \int_{M}X_{k}(f)d\mu=0, \quad \forall f\in C^\infty(M)
\ee
and that they satisfy the H\"ormander condition
\be\label{hyp_crochet}
\forall x\in M,\;\G_x=T_xM.
\ee
Let $\r\in \N$ be the smallest integer such that for any $x\in M$, $\G_x$ is generated by commutators of length at most $\r$.
For $k=1,\ldots, p$ and $x_0\in M$, denote
$\R\ni t\mapsto e^{tX_k}x_0$ the integral curve of $X_k$ starting from $x_0$ at $t=0$.\\

Let $h\in ]0,h_{0}]$ be a small parameter. Let us consider the following simple random walk 
 $x_{0},x_{1},...,x_{n},...$ on $M$, starting at $x_{0}\in M$: 
at step $n$, choose $j\in \{1,...,p\}$ at random and $t\in [-h,h]$ at random (uniform),
and set $x_{n+1}=e^{tX_j}x_n$. \\
\bigskip

Due to the condition $div(X_{j})=0$,
this random walk is reversible for the probability $\mu$ on $M$. It is easy to compute the Markov operator $T_{h}$ associated to this random walk:
 for any bounded and measurable function $f: M\rightarrow \R$ define
\be\label{def_op}
T_{k,h}f(x)=\frac 1{2h}\int_{-h}^hf(e^{tX_k}x)dt
\ee
Since the vector fields $X_k$ are divergence free,  for  any  $f,g$, we have
$$\int_MT_{k,h}f(x)g(x)d\mu=\int_Mf(x)T_{k,h}g(x)d\mu.$$
and the Markov operator associated to our random walk is
\be\label{def_op2}
T_hf(x)=\frac 1 p\sum_{k=1}^pT_{k,h}f(x)
\ee
One has $T_h(1)=1$, $\|T_h\|_{L^\infty\rightarrow L^\infty}=1$,
and $T_h$ can be uniquely extended as a bounded self-adjoint operator on $L^2=L^2(M,d\mu)$ 
such that $\|T_h\|_{L^2\rightarrow L^2}=1$.
In the following, we will denote $t_h(x,dy)$ the distribution kernel of $T_h$, and $t_h^n$ the kernel of $T_h^n$. Then, by construction, the probability for the walk starting at $x_{0}$
to be in a Borel set $A$ after $n$ step is equal to

$$ P(x_{n}\in A)=\int_{A} t_h^n(x_{0},dy)$$

The goal of this paper is to study the spectral theory of the operator $T_h$ and the convergence of $t_h^n(x_{0},dy)$ towards $\mu$ as $n$ tends to infinity.
Since $T_{h}$ is Markov and self adjoint, its spectrum is a subset of $[-1,1]$.
We shall denote by $g(h)$ the spectral gap of the  operator $T_{h}$. It is defined as the best constant
such that the following inequality holds true for all $u\in L^2$

\be\label{def_gap}
\Vert u \Vert^2_{L^2}-\<u, 1\>^2_{L^2}\leq {1\over g(h)}\<u-T_{h}u, u\>_{L^2}
\ee

The existence of a non zero spectral gap means that 
$1$ is a simple eigenvalue of $T_{h}$, and the distance between $1$ and the rest 
of the spectrum is equal to $g(h)$.
Our first result is the following

\begin{thm}\label{thm1} There exists $h_0>0$, $\delta_1,\delta_2>0$, $A >0$,  
and constants $C_i>0$ such that for  any $h\in ]0,h_0]$,
 the following holds true.\\

i) The spectrum of $T_{h}$ is a subset of $[-1+\delta_1,1]$, 
$1$ is a simple eigenvalue of $T_{h}$, 
and $Spec(T_{h})\cap [1-\delta_2,1]$ is discrete. 
Moreover, for any $0\leq \lambda \leq \delta_{2}h^{-2}$, the number of eigenvalues of $T_{h}$ in 
$ [1-h^2\lambda,1]$ (with multiplicity) is bounded by $C_1(1+\lambda)^{A}$.\\

ii) The spectral gap satisfies 

\be\label{gap3}
C_2 h^2 \leq g(h) \leq C_3 h^2
\ee
and the following estimate holds true for all integer $n$

\be\label{1.7}
sup_{x\in \Omega}\Vert t^n_{h}(x, dy)-\mu\Vert_{TV} \leq C_4e^{-n g(h) }
\ee
Here, for two probabilities on $M$, $\Vert \nu-\mu\Vert_{TV}=\sup_{A}\vert \nu(A)-\mu(A)\vert$
where the $\sup$ is over all Borel sets $A$, is the total variation distance between 
$\nu$ and $\mu$.
\end{thm}

We describe now the spectrum of $T_h$ near $1$. Let $\mathcal H^1(\X)$ be the Hilbert space
$$ \mathcal H^1(\X)=\{u \in L^2(M), \ \forall j=1,\ldots,p, \ X_{j}u\in L^2(M)\}$$
Let $\nu$ be the best constant such that the
following inequality holds true for all $u\in \mathcal H^1(\X)$ 
\be\label{nu1}
\Vert u \Vert^2_{L^2}-\<u, 1\>^2_{L^2}\leq {\mathcal E (u)\over \nu},
\ee
where 
\be
\mathcal E (u)= {1\over 6 }
\int_{M}\sum_{k=1}^p\vert X_{k} u\vert^2 d\mu
\ee
By the hypoelliptic theorem of H\"ormander (see \cite{Ho85}, Vol 3), one has 
$\mathcal H^1(\X)\subset H^s(M)$, for some $s >0$.
On the other hand, standard Taylor expansion in formula \eqref{def_op} show that for any fixed smooth function $g\in C^\infty( M)$, one has the following convergence in the space $C^\infty( M)$ 
\be\label{idel}
\lim_{h\rightarrow 0} {1-T_{h} \over h^2} g =L(g),
\ee
where  the operator $L=-{1\over 6p}\sum_{k} X_{k}^2$ is the positive Laplacian associated to the Dirichlet form $\mathcal E (u)$.
It has a compact resolvant and spectrum  $\nu_0=0<\nu_1=\nu <\nu_2 < ...$.
Let $m_j$ be the multiplicity of $\nu_j$. One has $m_0=1$ since $Ker(L)$ is spaned by the
constant function $1$ thanks to the Chow theorem (\cite{Chow39}). In fact, for any $x,y\in M$
there exists  a continuous  curve connecting $x$ to $y$  which is
a finite union of pieces of trajectory of one of the fields $X_{j}$.\\

\begin{thm}\label{thm2}
One has
\be\label{1.6}
lim_{h\rightarrow 0} h^{-2}g(h)=\nu
\ee
Moreover, for any $R>0$ and $\varepsilon>0$ such that the intervals $[\nu_j-\varepsilon, \nu_j+\varepsilon]$ are disjoint for $\nu_{j}\leq R$,   there exists $h_1>0$ such that for all $h\in ]0,h_1]$
\be\label{1.5}
Spec({1-T_{h}\over h^2})\cap ]0,R] \subset \cup_{j\geq 1} [\nu_j-\varepsilon, \nu_j+\varepsilon]
\ee
 and the number of eigenvalues of ${1-T_h\over h^2}$ with multiplicities, in the interval 
 $[\nu_j-\varepsilon, \nu_j+\varepsilon]$,  is equal to $m_j$. 
\end{thm}

The paper is organized as follows:\\

In section \ref{sec1}, we  recall some basic facts on nilpotent Lie groups, and we recall the Goodman version (see \cite{Go78}) of one of the main results of the Rotschild-Stein paper.\\
  
In section \ref{sec2}, the main result is the proposition \ref{prop:minor_iter} which
gives a lower bound on a suitable power $T_{h}^P$ of $T_{h}$. This in particular allows to get
a first crude but fundamental bound on the $L^\infty$ norms of eigenfunctions of 
$T_{h}$ associated to eigenvalues close to $1$. \\

Section \ref{sec4} is devoted to the study of the Dirichlet form associated
to our random walk. The fundamental result of this section is  proposition \ref{prop4.1}.
It allows to separate clearly the spectral theory of $T_{h}$ in low and high frequencies with
respect to the parameter $h$. Many tools in this section are essentially an adaptation 
to the semi-classical setting of the ideas contained in the  Rotschild-Stein paper\cite{RoSt76}.\\

Section \ref{TV} is devoted to the proof of theorems \ref{thm1} and \ref{thm2}.
With propositions \ref{prop:minor_iter} and \ref{prop4.1} in hands, the  proof follows 
the strategy of \cite{DiLeMi11} and \cite{DiLeMi12}, but with some differences and simplifications.
This section  contains also a paragraph on the Fourier analysis associated to $T_{h}$
that will be useful in \ref{sec6}. \\

Section \ref{sec6} is devoted to the proof of
the convergence when $h\rightarrow 0$ of our Markov chain to the hypoelliptic diffusion on the manifold
$M$
associated to the generator $L={-1\over 6p}\sum_{k}X_{k}^2$ . This is probably a well known result for specialists, but we have not succeed 
to find a precise reference. Since this convergence follows as a simple
byproduct of our estimates, we decide to include it in the paper. \\

Finally, the appendix contains two lemmas. Lemma \ref{lemA1} shows how to deduce from
proposition \ref{prop4.1} a Weyl type estimate on the eigenvalues of $T_{h}$
in a neighborhood of $1$. Lemma \ref{lemcoho} is an elementary cohomological lemma
on the Schwartz space of the nilpotent Lie algebra $\mathcal N$.\\

{\bf Acknowledgement:} We thank Dominique Bakry who has motivated us to study this problem. 

\newpage
\section{The lifted operator to a nilpotent Lie algebra}\label{sec1}

We will use the notation $\N_q=\{1,\ldots,q\}$.
For any  family of vector fields $Z_1,\ldots, Z_p$ and any multi-index $\alpha=(\alpha_1,\ldots,\alpha_k)\in\N_p^k$ 
denote $|\alpha|=k$ the length of $\alpha$ and let 
\be\label{2.1}
Z^\alpha=H_{\alpha}(Z_{1},...,Z_{p})=[Z_{\alpha_1},[Z_{\alpha_2},\ldots[Z_{\alpha_{k-1}},Z_{\alpha_k}]\ldots]
\ee

Let $\Y_1,\ldots,\Y_p$ be a system of generators of the free lie algebra with $p$ generators $\F$ and let $\A^\infty$ be a set of multi-indexes such that 
$(\Y^\alpha)_{\alpha\in \A^\infty}$ is a basis of $\F$.

Introduce
$\NN$ the free up to step $\r$ nilpotent Lie algebra generated by $p$ elements $Y_1,\ldots, Y_p$, and let $N$ be the corresponding simply connected Lie group.
We have the decomposition
\be
\NN=\NN_1\oplus\ldots\oplus\NN_\r
\ee
where $\NN_1$ is generated by $Y_1,\ldots Y_p$ and  
$\NN_{j}$ is spanned by the commutators $Y^\alpha=H_{\alpha}(Y_{1},...,Y_{p})$
with $\vert \alpha \vert=j$ for $2\leq j\leq \r$. Let 
$\A=\{\alpha\in \A^\infty,\, |\alpha|\leq \r\}$ and $\A_r=\{\alpha\in \A,\;|\alpha|=r\}$.
The family
 $(Y^\alpha)_{\alpha\in \A}$ is a  basis for $\NN$ and for any $r\in\N_\r$, $\{Y^\alpha,\; \alpha\in \A_r\}$ is a basis of $\NN_r$.
We denote by  $D=\sharp\A$ the dimension of $\NN$.
 The action of $\mathbb R_{+}$ on $\NN$ is given by
$$t.(v_{1},v_{2},...,v_{r})=(tv_{1},t^2v_{2},...,t^\r v_{\r})$$
An homogeneous norm  $\bg v \bg $ which is smooth in $\NN \setminus o_{\NN}$
is given by 
$$\bg v \bg=(\sum_{j}\vert v_{j}\vert^{(2\r!)/j})^{1/(2\r!)}$$
 where
$\vert v_{j}\vert$ is an euclidian norm on $\NN_{j}$, and   
$$ Q=\sum j \dim(\NN_{j})$$
is the quasi homogeneous dimension of $\NN$.
We will identify the Lie agebra $\NN$ with the Lie group $N$ by the exponential map, i.e  
the product law $a.b$  on $\NN$ is given  by $exp(a.b)=exp(a)exp(b)$.
In particular, one has with this identification $a^{-1}=-a$ for all $a\in \NN$.
To avoid notational confusion, we will use sometime the notation $e=o_{\NN}$, so that
$a.e=e.a=a$ for all $a\in \NN$. For $Y\in T_{e}\NN\simeq \NN$, we denote by 
$\tilde Y$ the left invariant vector field on $\NN$ such that $\tilde Y(o_{\NN})=Y$, i.e

$$\tilde Y(f)(x)= {d\over ds}(f(x.sY)\vert_{s=0}$$
The right invariant vector field on $\NN$ such that $Z(o_{\NN})=Y$ is defined by
$$Z(f)(x)= {d\over ds}(f(sY.x)\vert_{s=0}$$
Here, $sY$ is the usual product of the vector $Y\in \NN$ by the scalar $s\in \R$.
For $a\in \NN$, let $\tau_{a}$ be the diffeomorphism of $\NN$ defined by $\tau_{a}(u)=a.u$.
One has 
$$\tilde Y(a)=d\tau_{a}(e)(Y)$$ 
\begin{exa} The standard 3d-Heisenberg group is $\NN=\R^2_{x,y}\oplus \R_{t}\simeq \R^3$,
with the product law
$$(x,y,t).(x',y',t')=(x+x',y+y',t+t'+xy'-yx')$$
and the left invariant vector fields associated respectively to the vectors
$(1,0,0), (0,1,0)$ and $(0,0,1)$ are in that case
$$\tilde Y_{1}= {\partial\over \partial x}-y{\partial\over \partial t},\quad 
\tilde Y_{2}={\partial\over \partial y}+x{\partial\over \partial t}, \quad \text{and} \quad
{\partial\over \partial t}={1\over 2}[\tilde Y_{1},\tilde Y_{2}]$$
\end{exa}

\noindent

In general, for $x=(x_{1},..., x_{\r})$  and 
$y=(y_{1},..., y_{\r})$, $x_{j}, y_{j}\in \NN_{j}$, the product law is given by
\be\label{3.10}
\ba
&(x_{1},..., x_{\r}).(y_{1},..., y_{\r})=(z_{1},..., z_{\r})\\
& z_{j}=x_{j}+y_{j}+P_{j}(x_{<j},y_{<j})
\ea\ee
with the notation $x_{<j}=(x_{1},...,x_{j-1})$, and where $P_{j}$ is a polynomial
of degree $j$ with respect to the homogeneity on $\NN$, i.e
$$P_{j}((t.x)_{<j},(t.y)_{<j})=t^jP_{j}(x_{<j},y_{<j})$$
which is compatible with the identity $t.(x.y)=(t.x).(t.y)$. \\

Let $\lambda: \NN\rightarrow  \G$ be the unique linear map such that for any $\alpha\in \A$, $\lambda(Y^\alpha)=X^\alpha$. Then $\lambda$ is a Lie homomorphism ``up to step $\r$'':
\be
\lambda([Y^\alpha,Y^\beta])=[X^\alpha,X^\beta]
\ee
for any multi-indexes $\alpha,\beta$ such that $|\alpha|+|\beta|\leq \r$.

Let $x_0\in M$. There exists a subset $\A_{x_0}\subset\A$ such that 
$(X^\alpha(x))_{\alpha\in \A_{x_0}}$ is a basis of $T_{x}M$ for any $x$ close to $x_{0}$.
Therefore,  there exists a neighborhood $\Omega_0$ of the origin $o_\NN$ in $\NN$ 
and a neighborhood $V_0$ of $x_{0}$ in $M$ such that the map $\Lambda$
$$\Lambda: u=\sum_{\alpha\in \A}u_{\alpha}Y^\alpha\in\Omega_0\mapsto e^{\lambda(u)}x_0=e^{\sum_{\alpha\in \A}u_{\alpha}X^\alpha}x_{0}$$
 is a  submersion from $\Omega_0$ onto $V_0$,
and the map 
$W_{x_0}:\mathcal{C}^\infty(V_0)\rightarrow \mathcal{C}^\infty(\Omega_0)$
defined by $W_{x_0}f(u)=f(e^{\lambda(u)}x_0)$
is injective. Since $\Lambda$ is a submersion,  
there exists a system of coordinates $\theta:\R^m\times\R^n\rightarrow\NN$ 
defined near $o_\NN$, where $m+n=D$, such that $\Lambda\theta:\R^m\rightarrow M$ is a system 
of coordinates near $x_{0}$ and in these coordinates one has $\Lambda (x,y)=x$. We thus may assume
that in these coordinates one has  $\Omega_{0}=V_{0}\times U_{0}$ where $U_{0}$ is a neighborhood of $0\in \R^n$.\\
\begin{exa}\label{exa1}
Take for example the two vectors fields in $\R^2$, $X_{1}=\partial_{x}$, 
$X_{2}=x\partial_{y}$ so that $[X_{1},X_{2}]=\partial_{y}$. Then on take for 
$\mathcal N$ the $3d$-Heisenberg group, and the map $\lambda$ is given by, with $T=2\partial_{t}=[Y_{1},Y_{2}]$
$$\lambda(u_{1}Y_{1}+u_{2}Y_{2}+u_{3}T)=
u_{1}X_{1}+u_{2}X_{2}+u_{3}[X_{1},X_{2}]=u_{1}\partial_{x}+ (u_{3}+u_{2}x)\partial_{y}$$
Thus we get
\be\label{3.10bis}
e^{\lambda(u)}(x,y)=(x+u_{1}, y+u_{3}+u_{2}x+{1\over 2}u_{1}u_{2})
\ee
Let $I_{h}=\{\vert u_{1}\vert < h,\vert u_{2}\vert < h,
\vert u_{3}\vert < h^2\}$. One has $Vol(I_{\epsilon,h})=8h^4$.
Observe on this example that the set 
$\tilde B_{h,(x,y)}=\{e^{\lambda(u)}(x,y), u\in I_{h}\}$, when $(x,y)$ is fixed
 and $h$ small, has volume of order:\\
  $h^2$ when $x\not=0$, and $h^3$ when $x=0$.

\end{exa}

Let us now recall the notion of order of a vector field used in \cite{RoSt76} and \cite{Go78}. Denote $\{\delta_t\}_{t>0}$ the one parameter group of dilating automorphisms on $\NN$:
$$\delta_t Y^\alpha=t^{|\alpha|}Y^\alpha.$$
Let $\Omega$ be a compact neighborhood of $o_{\NN}$ in $\NN$. For any $m\in\N$, let 
$$C^\infty_m=\{u\in C^\infty( \Omega,\R),\, f(u)=\O(\bg u\bg^m)\}.$$
We have the filtration $C^\infty(\Omega)=C^\infty_0\supseteq C^\infty_1\supseteq\ldots$, and $C^\infty_m . C^\infty_n\subseteq C^\infty_{m+n}$.
Let $T: C^\infty(\Omega)\rightarrow C^\infty(\Omega)$. We say that $T$ is of order less than $k$ at $0$, if $T(C^\infty_m)\subseteq C^\infty_{m-k}$ for all integers $m\geq 0$.
If $\partial_\alpha$ denotes the differentiation in the direction $Y^\alpha$, then a vector field $T=\sum_\alpha\varphi_\alpha \partial_\alpha$ is of order $\leq k$ iff 
$\varphi_\alpha\in C^\infty_{|\alpha|-k}$ for all $\alpha$, with the convention 
$C^\infty_m=C^\infty_0$ for $m\leq 0$.

The following result is the Goodman version of one of the results of the article
\cite{RoSt76} by  L.Rothschild and E.Stein.
\begin{thm}\label{goodman}
Decreasing $\Omega_{0}$ if necessary, there exists $C^\infty$ vector fields $Z_1,\ldots, Z_p$ on $\Omega_{0}$ such that for any $\alpha\in \A$ we have
\begin{enumerate}
 \item[i)] $Z^\alpha W_{x_{0}}=W_{x_{0}}X^\alpha$
\item[ii)] $Z^\alpha=\tY^\alpha+R_\alpha$, where $R_\alpha$ is a vector field of order $\leq|\alpha|-1$ at $0$.
\item[iii)] $Z^\alpha=H_{\alpha}(Z_{1},...,Z_{p})$ for all ${\alpha\in\A}$.
\end{enumerate}

\end{thm}

Observe that in the previous coordinate system $(x,y)$ on $\Omega_{0}$, one can write
for $\alpha \in \mathcal A$
\be\label{2.0}
 X^\alpha=\sum_{j}a_{\alpha,j}(x){\partial \over \partial x_{j}}, \quad
Z^\alpha=\sum_{j}a_{\alpha,j}(x){\partial \over \partial x_{j}}+ \sum_{l}b_{\alpha,l}(x,y)
{\partial \over \partial y_{l}}
\ee

As an obvious  consequence of this theorem, we have the following, with $W=W_{x_{0}}$,
and  $\tilde\lambda (u)=\sum_{\alpha\in A}u_{\alpha}Z^\alpha$.
\begin{prop}\label{prop2.2}
 Let $f\in \mathcal{C}^0(V_0)$ and let $\omega_0\subset\subset\Omega_0$ be a neighborhood of $o_{\NN}$. Then, there exists $r_0>0$ such that 
for all $\bg u \bg\leq r_{0}$,  and  $v\in \omega_0$, we have
\be
(Wf)(e^{\tilde\lambda (u)}v)=W(f_{u})(v)
\ee
where the function $f_{u}$ 
 is defined near $x_{0}$ by $f_{u}(x)=f(e^{\lambda(u)}x)$.
\end{prop}

Using this proposition, we can easily compute the action of $W$ on the operator $T_h$ acting on functions with support close to $x_{0}$.  We get immediatly
\be\label{act_W_T1}
WT_{h}=\tT_{h}W, \quad \tT_{h}={1\over p}\sum_{k=1}^p \tT_{k,h}
\ee
where for  $u\in\NN$ small.
\be\label{def_tT}
\tT_{k,h}g(u)=\frac 1 {2h}\int_{-h}^hg(e^{tZ_k}u)dt
\ee
  Using the notation $T^\alpha=T_{\alpha_k,h}\ldots T_{\alpha_1,h}$ for any multi-index $\alpha=(\alpha_1,\ldots,\alpha_k)$ we get for any $u\in \NN$ close to 
 $o_{\NN}$ such that
 $\Lambda(u)=x$ 
\be\label{act_W_T}
T^\alpha f(x)=W(T^\alpha f)(u)=\frac 1 {(2h)^k}\int_{[-h,h]^k}(Wf)(e^{t_1Z_{\alpha_1}}\ldots e^{t_kZ_{\alpha_k}}u)dt_1\ldots dt_k
\ee

\section{Rough bounds on eigenfunctions}\label{sec2}

Let us recall from section \ref{sec1}, that for $u=\sum_{\alpha\in\mathcal A}u_{\alpha}Y^\alpha\in \NN$, the vector field $\lambda(u)$ on $M$ is defined by  $\lambda(u)=\sum_{\alpha\in\mathcal A}u_{\alpha}X^\alpha$. Let $\epsilon >0$
 and 
$I_{\epsilon,h}$ be the neighborhood of $o_{\NN}$ in $\NN$ defined by
$$I_{\epsilon,h}=\{u=\sum_{\alpha\in\mathcal A}u_{\alpha}Y^\alpha, \quad u_{\alpha}\in ]-\epsilon h^{|\alpha|}, \epsilon h^{|\alpha|}[ \ \} $$ 

For any $x\in M$ we define a positive measure $S^\epsilon_h(x,dy)$ on $M$
by the formula

\be\label{def_Sj}
\forall f\in C^0(M), \quad \int f(y)S^\epsilon_{h}(x,dy)=h^{-Q}\int_{u\in I_{\epsilon,h}}f(e^{\lambda (u)}x)\ du
\ee
where $du=\Pi_{\alpha}du_{\alpha}$ is the left and right invariant Haar measure on $\NN$. 
Let us introduce the numerical sequence $(b_n)_{n\in\N^*}$ defined by $b_1=1$ and $b_{n+1}=2b_n+2$, so that for all $n\in\N^*$, we have $b_n=3.2^{n-1}-2$.
For all $r=1,\ldots,\r$ denote $a_r=\sharp\A_r=\dim\NN_r$, and let $P=\sum_{r=1}^\r a_rb_r$.

\begin{prop}\label{prop:minor_iter}
 There exists  $\epsilon >0$, $c>0$ and $h_0>0$ such that 
for all $h\in]0,h_0]$,  $x\in M$
\be\label{iter_noyau}
t_h^P(x,dy)=\rho_h(x,dy)+cS^\epsilon_{h}(x,dy)
\ee
where $\rho_h(x,dy)$ is a non-negative Borel measure on $M$ for all $x\in M$.
\end{prop}
\begin{remk}
As in \cite{DiLeMi11}, one can deduce from proposition \ref{prop:minor_iter} that the inequality \eqref{iter_noyau} holds true for $t_h^N(x,dy)$ as soon as $N\geq P$,
eventually with different constants $\epsilon >0$, $c>0$ and $h_0>0$  depending on $N$.
\end{remk}

Before proving this proposition, let us give two simple but fundamental corollaries. 
Like in \cite{DiLeMi11}, these two corollaries will play a key role in the proof of theorems 
\ref{thm1} and \ref{thm2}.
Here, we use the same notation for a bounded measurable family in $x$ of
non negative  Borel measure $k(x,dy)$ and the corresponding operator $f\mapsto K(f)(x)=\int f(y)k(x,dy)$ acting on $L^\infty$.

\begin{cor}\label{cor:borne_spec_ess}
 There exists $h_0>0$ and $\gamma<1$ such that for all $h\in]0,h_0]$ and all $x\in M$
\be
\|\rho_h(x,dy)\|_{L^\infty\rightarrow L^\infty}\leq \gamma <1
\ee
\end{cor}

\bp
By definition, the non-negative measure $\rho_h$ is given by $\rho_h(x,dy)=t_h^P(x,dy)-cS^{\epsilon}_h(x,dy)$. Therefore
\be\label{borne_inf_mu}
|\int_Mf(x)d\rho_h(x,dy)|\leq \|f\|_{L^\infty}\int_Md\rho_h(x,dy)
\leq \|f\|_{L^\infty}(1-c\inf_{x\in M}\int_MS^{\epsilon}_h(x,dy))
\ee
since $t_h^P(x,dy)$ is a Markov kernel.
From \eqref{def_Sj}, one has $\int_MS^{\epsilon}_h(x,dy)=h^{-Q}\operatorname{meas}(I_{\epsilon,h})=(2\epsilon)^D$.
Combined with \eqref{borne_inf_mu}, this implies the result.
\ep

\begin{cor}\label{cor:borne_L2_Linf}
Let $a\in]\gamma^{\frac 1 P},1]$ be fixed. There exists $C=C_a>0$ such that for any $\lambda\in[a,1]$ and any $f\in L^2(M,d\mu)$ we have
\be\label{L2Linf}
T_hf=\lambda f\Longrightarrow \|f\|_{L^\infty}\leq Ch^{-\frac Q 2}\|f\|_{L^2}
\ee
\end{cor}

\bp
Suppose $T_hf=\lambda f$, then $T_h^Pf=\lambda^Pf$. Hence, $S^\epsilon_hf=\lambda^Pf-\rho_h(f)$ and then
\be\label{borne_L2_Linf_1}
\|S^\epsilon_hf\|_{L^\infty}\geq \lambda^P\|f\|_{L^\infty}-\gamma\|f\|_{L^\infty}\geq c_a\|f\|_{L^\infty}
\ee
with $c_a=a^P-\gamma$.
On the other hand, since $u\mapsto e^{\lambda(u)}x$ is a submersion from a
neighborhood of $o_{\NN}\in \NN$ onto a neighborhood of $x\in M$, we get by Cauchy-Schwarz inequality

\be\label{borne_L2_Linf_2}
 \vert S_h^\epsilon f(x)\vert \leq h^{-Q}\operatorname{meas}(I_{\epsilon,h})^{1/2}
 (\int_{u \in I_{\epsilon,h}}|f(e^{\lambda(u)}x)|^2\ du)^{1/2} \leq Ch^{-Q/2}\Vert f\Vert_{L^2(M)}
\ee

Putting together \eqref{borne_L2_Linf_1} and \eqref{borne_L2_Linf_2}, we obtain the anounced result.
\ep

Le us now prove Proposition \ref{prop:minor_iter}.
We have to show that there exists $c,\epsilon>0$ independent of $h$ small, such that for any non negative continous function $f$ on $M$, one has
$T_h^Pf(x)\geq cS^\epsilon_hf(x)$. 
Since $M$ is compact and the operator $T_h$ moves supports of functions at distance at most $h$,  we can assume without loss of generality that $f$ is  supported near some point $x_0\in M$ where we can apply the results of section \ref{sec1}. 
Recall $\tilde\lambda (u)=\sum_{\alpha\in A}u_{\alpha}Z^\alpha$. From proposition \ref{prop2.2} one has
$f(e^{\lambda (u)}x)=W(f)(e^{\tilde\lambda (u)}w)$ for any $w$ close to $o_{\NN}$
such that $\Lambda(w)=x$. Using also \eqref{act_W_T1},
we are thus reduce  to prove that there exists
$c,\epsilon>0$ independent of $h$ small such that for any  non negative continuous function $g$ on $\NN$
supported near $o_{\NN}$, one has 
\be\label{2.20}
\tT_{h}^P g (w)\geq c h^{-Q}\int_{u\in I_{\epsilon,h}}g(e^{\tilde\lambda (u)}w)\ du
\ee

For any non-commutative sequence $(A_k)$ of operators, we denote $\Pi_{k=1}^KA_k=A_K\ldots A_1$ (i.e $A_{1}$ is the first operator acting). Endowing $\A_r$ with the lexicographical order, we can write $\A_r=\{\alpha_1<\ldots<\alpha_{a_r}\}$ and for any non-commutative sequence $(B_\alpha)$ indexed by $\A$, we define
$\Pi_{\alpha\in\A_r}B_\alpha=\Pi_{j=1}^{a_r}B_{\alpha_j}$ and 
$\Pi_{\alpha\in \A}B_\alpha=\Pi_{r=1}^\r\Pi_{\alpha\in\A_r}B_\alpha$.\\

Let $\alpha=(\alpha_{1},...,\alpha_{k})\in \mathbb N_{p}^k$, 
and $t=(t_{1},...,t_{k})\in \R^k$ close to $0$. One defines by induction on $\vert\alpha\vert$ a smooth diffeomorphism
$\phi_{\alpha}(t)$ of $\NN$ near $o_{\NN}$, with $\phi_{\alpha}(0)=Id$, by the following formulas:\\
If $\vert\alpha\vert=1$ and $\alpha=j\in \{1,...,p\}$, set
$\phi_{\alpha}(t)(w)=e^{tZ_{j}}w$.
If $\vert\alpha\vert=k\geq 2$, set $\alpha=(j,\beta)$, with $\beta\in \mathbb N_{p}^{k-1}$
and $t=(t_{1},t')$ with $t'\in \R^{k-1}$ and set
\be\label{2.10bis}
\phi_{\alpha}(t)=\phi_{\beta}^{-1}(t')e^{-t_{1}Z_{j}}\phi_{\beta}(t')e^{t_{1}Z_{j}}
\ee
Observe that $\phi_{\alpha}(t)=Id$ if one of the $t_{j}$ is equal to $0$. The map $(t,w)\mapsto \phi_{\alpha}(t)(w)$ is smooth, and one has in local coordinates
on $\NN$ and for $t$ close to $0$
\be\label{2.10}
\phi_{\alpha}(t)(w)= w+ (\Pi_{1\leq l\leq \vert\alpha\vert}\ t_{l})\ Z^\alpha(w)+ r_{\alpha}(t,w)
\ee
with $r_{\alpha}(t,w)\in (\Pi_{1\leq l\leq \vert\alpha\vert}\ t_{l})O(\vert t \vert)$. From \eqref{2.10bis}, one get
easily by induction on $k$ the following lemma.

\begin{lem}\label{lemfi}
For  $2\leq k \leq \r$, there exists maps:
$$\epsilon_k:\,\{1,\ldots,b_k\}\rightarrow\{\pm 1\},\;\;\ell_k:\,\{1,\ldots,b_k\}\rightarrow \{1,\ldots,k\},\;j_k:\,\{1,\ldots,b_k\}\rightarrow\{1,\ldots,p\}$$ 
such that $\epsilon_k(1)=1$, $\epsilon_k(b_k/2)=-1$, $\ell_k(1)=1$, $\ell_k(b_k/2)=1$, $\sharp \ell_{k}^{-1}(j)=2^j$ for $j\leq k-1$, $\sharp \ell_{k}^{-1}(k)=2^{k-1}$,
$j_{k}(m)=\alpha_{\ell_k(m)}$, and such that for all  $t=(t_{1},...,t_{k})$ one has 
\be
\phi_\alpha(t)=\prod_{m=1}^{b_k}e^{\epsilon_k(m)t_{\ell_k(m)}Z_{j_k(m)}}
\ee
\end{lem}

Since $g$ is non negative, one has 
\be
\tT_h^Pg(w)\geq \frac 1{p^P}\prod_{\alpha\in\A}\prod_{k=1}^{b_{|\alpha|}}T_{j_{\vert\alpha\vert}(k),h}g(w)\\
\ee
Therefore,  we are reduced to prove that there exists $\epsilon,c>0$ independent
of $h$ small and $w$ near $o_{\NN}$ such that the following inequality holds true.
\be\label{2.20bis}
h^{-P}\int_{[-h,h]^P}g\Big(\prod_{\alpha\in\A}\prod_{k=1}^{b_{|\alpha|}}
e^{t_{\vert \alpha\vert,k}Z_{j_{\vert\alpha\vert}(k)}}w\Big)dt
\geq c h^{-Q}\int_{z\in I_{\epsilon,h}}g(e^{\tilde\lambda (z)}w)\ dz
\ee

Let $\Phi_{w}:\R^{P}\longrightarrow\NN$ be the smooth map defined for $s=(s_{\alpha,k})_{\alpha\in\A,k=1,\ldots,b_{|\alpha|}}\in \R^P$ by the formula
\be\label{def_psi}
\Phi_{w}(s)= \Big(\prod_{r=1}^\r\prod_{\alpha\in\A_{r}}\prod_{k=1}^{b_r}e^{s_{\alpha,k}
Z_{j_{\vert\alpha\vert}(k)}}\Big)w
\ee
Since $(Z^\beta(w))_{\beta\in \A}$ is a basis of $T_{w}\NN$, 
$u=(u_{\beta})_{\beta\in\A}\mapsto e^{\sum_{\beta\in\A}u_{\beta}Z^\beta}w$
is a local coordinate system centered at $w\in \NN$, and therefore, there exists
smooth functions $U_{\beta,w}(s)$ such that
\be\label{2.11}
\Phi_{w}(s)=e^{\sum_{\beta\in\A}U_{\beta,w}(s)Z^\beta}w
\ee 
Moreover, it follows easily from the Campbell-Hausdorff formula, that one has $U_{\beta,w}(s)\in O(s^{\vert\beta\vert})$
near $s=0$. Let now $\kappa:\R^{Q}\longrightarrow \R^{P}$ the map defined by

\be\label{def_kappa}
(t_{\alpha,l})_{\alpha\in\A,l\in\N_{|\alpha|}}\mapsto 
(\epsilon_\alpha(k)t_{\alpha,\ell_{\vert\alpha\vert}(k)})_{\alpha\in\A,k=1,\ldots,b_{|\alpha|}}
\ee
Then, from lemma \ref{lemfi} we have the following identity for any $t=(t_{\alpha})_{\alpha\in \A}\in\R^Q$.
\be\label{compos_psi_kappa}
\Phi_{w}\circ\kappa(t)=\Pi_{\alpha\in\A}\phi_{\alpha}(t_{\alpha})w
\ee
From \eqref{2.10} and Campbell-Hausdorff formula, one gets
\be\label{2.15}
\ba
&\Pi_{\alpha\in\A}\phi_{\alpha}(t_{\alpha})w=e^{\sum_{\beta\in \A}f_{\beta}(t)Z^\beta}w\\
&f_{\beta}(t)=\Pi_{1\leq l \leq \vert \beta\vert}t_{\beta,l}+ g_{\beta}((t_{\gamma})_{\vert\gamma\vert<\vert\beta\vert})+ r_{\beta}(t)
\ea\ee
with $g_{\beta}$ an homogeneous polynomial of degree $\vert\beta\vert$ depending 
only on $(t_{\gamma})_{\vert\gamma\vert<\vert\beta\vert}$ and 
$r_{\beta}(t)\in O(\vert t \vert^{\vert\beta\vert+1})$. Let $ \delta\in]\frac 1 2, 1[$ and define $\xi=(\xi_{\alpha,k})_{\alpha\in\A,k\in\N_{|\alpha|}}\in\R^Q$ by $\xi_{\alpha,1}=0$ and $\xi_{\alpha,k}=\delta h$ for  $k= 2,\ldots,\vert \alpha\vert$.
Let  $\zeta:\R^D\longrightarrow\R^Q$ be the map defined by the formula
\be
\ba
&s=(s_\alpha)_{\alpha\in\A}\mapsto (\zeta_{\alpha,k}(s))_{\alpha\in\A,k\in\N_{|\alpha|}}\\
&\zeta_{\alpha,1}(s)=s_\alpha, \quad \text{and}\quad  \zeta_{\alpha,k}(s)=0 \quad\forall k\geq 2
\ea\ee
and let $\sigma:\R^{P-D}\longrightarrow \R^P$ be the map defined by the formula
\be
\ba
& v=(v_{\alpha,k})_{\alpha\in\A,k=2,\ldots,b_{|\alpha|}}\mapsto (\sigma_{\alpha,k}(v))_{\alpha\in\A,k=1,\ldots,b_{|\alpha|}}\\
&\sigma_{\alpha,1}(v)=0,  \quad \text{and}\quad \sigma_{\alpha,k}(v)=v_{\alpha,k}
\quad \forall k\neq 1
\ea\ee
Set $\hat\kappa_{\xi}(u,v)=\kappa(\zeta(u)+\xi)+\sigma(v)$, and let $\Psi_{w}:\R^D\times\R^{P-D}\rightarrow\NN$ be defined by
\be\label{def_Psi}
\Psi_{w}(u,v)=\Phi_{w}(\hat\kappa_{\xi}(u,v))
\ee
Then, it follows from \eqref{2.11} that there exists smooth maps 
$\hat\varphi_{\alpha,w}(u,v)$ such  that 
\be\label{log_Psi}
\Psi_{w}(u,v)=e^{\sum_{\alpha\in\A}\hat\varphi_{\alpha,w}(u,v)Z^\alpha}w
\ee
From \eqref{compos_psi_kappa} one has 
$$\Psi_{w}(u,0)=\Phi_{w}(\kappa(\zeta(u)+\xi))=
\Pi_{\alpha\in\A}\phi_{\alpha}(u_{\alpha},\delta h,...,\delta h)w$$ 
and therefore from  \eqref{2.15} we get, since $\hat\kappa_{\xi}(u,v)$ is linear in $\xi,u,v$

\be\label{eq:def_hatphi}
\hat\varphi_{\alpha,w}(u,v)=u_\alpha (\delta h)^{|\alpha|-1}+ 
g_{\alpha,w}((u_{\gamma})_{\vert\gamma\vert<\vert\alpha\vert},\delta h) +
p_{\alpha,w}(u,\delta h, v)+q_{\alpha,w}(u,\delta h,v)
\ee
where $g_{\alpha,w}(u,s)$ is a homogenous polynomial of degre $|\alpha|$ depending only on 
$u_{\gamma}$ for $\vert\gamma\vert<\vert\alpha\vert$,
$p_{\alpha,w}(u,s, v)$ is a homogenous polynomial of degre $|\alpha|$ in $(u,s,v)$
such that $p_{\alpha,w}(u,s,0)=0$, and $q_{\alpha,w}(u,s,v)\in O((u,s,v)^{1+\vert\alpha\vert})$
near $(u,s,v)=(0,0,0)$. Moreover, from $\phi_{\alpha}(0,\delta h,...,\delta h)=Id$, one get $g_{\alpha,w}(0,s)=0$ and also $q_{\alpha,w}(0,s,0)=0$. Observe that $w$ is just a smooth parameter in the above constructions. 
Thus, we will remove the dependance in $w$ in what follows. Define now
\be
\begin{split}
 \mathfrak Q:\R^P=\R^D\times\R^{P-D}&\longrightarrow \R^P\\
(u,v)=((u_\alpha)_{\alpha\in\A},(v_{\alpha,k})_{\alpha\in\A,k=2,\ldots,b_{|\alpha|}})&\mapsto((\hat\varphi_{\alpha}(u,v))_{\alpha\in\A}, v)
\end{split}
\ee
and for $\eta,\epsilon>0$ let $$\Delta_{\epsilon,\eta}=\{(u,v)=((u_\alpha)_{\alpha\in\A},(v_{\alpha,k})_{\alpha\in\A,k=2,\ldots,b_{\vert\alpha\vert}})\in\R^P,\;\vert u_\alpha\vert<\epsilon h,\text{ and }
\vert v_{\alpha,k}\vert <\eta h\text{ for all }\alpha, k\}$$

\begin{lem}\label{lem:diffeo_theta} Let $\delta\in]\frac 1 2 , 1[$  be fixed.
 There exists $0<\eta<<\epsilon<1/2$  and $h_0>0$ such that the restriction $\mathfrak Q_{\epsilon,\eta}$ of $\mathfrak Q$ to $\Delta_{\epsilon,\eta}$ enjoys the following:
\begin{enumerate}
 \item there exists $U_{\epsilon,\eta}$ open neighbourhood of $0\in\R^P$ such that $\mathfrak Q_{\epsilon,\eta}:\Delta_{\epsilon,\eta}\rightarrow U_{\epsilon,\eta}$ is a $C^\infty$ diffeomorphism.
\item there exists some constant $C>0$ such that for all $h\in]0,h_0]$ and all $(u,v)\in \Delta_{\epsilon,\eta}$
$$h^{Q-D}/C\leq \operatorname{J\mathfrak Q}_{\epsilon,\eta}(u,v):=|\det(D_{(u,v)}\mathfrak Q_{\epsilon,\eta})|\leq Ch^{Q-D}$$
\item there exists $M\geq 1$ such that for all $h\in]0,h_0]$ the set $U_{\epsilon,\eta}$ contains $I_{\epsilon/M,h}\times ]-\eta h,\eta h[^{P-D}$
where
$I_{\epsilon/M,h}=\prod_{\alpha\in \A}]-\epsilon h^{|\alpha|}/M,\epsilon h^{|\alpha|}/M[.$
\end{enumerate}
 
\end{lem}

\bp The proof is just a scaling argument. Set $u_{\alpha}=h\tilde u_{\alpha}$, 
$v_{\alpha,k}=h\tilde v_{\alpha,k}$ and $\hat\varphi_{\alpha}=h^{\vert\alpha\vert}z_{\alpha}$.
Then the map $\mathfrak Q$ becomes after scaling
$\tilde{\mathfrak Q}: (\tilde u,\tilde v)\mapsto (z,\tilde v)$, and from \eqref{eq:def_hatphi} one has 
$$ z_{\alpha}=\tilde u_\alpha \delta^{|\alpha|-1}+
g_{\alpha}((\tilde u_{\gamma})_{\vert\gamma\vert<\vert\alpha\vert},\delta )+
p_{\alpha}(\tilde u,\delta , \tilde v)+h\tilde q_{\alpha}(\tilde u,\delta,\tilde v,h)$$
$p_{\alpha}(\tilde u,\delta ,0)=0$,  $\tilde q_{\alpha}(\tilde u,\delta,\tilde v,h)$
is smooth and vanish at order $\vert\alpha\vert+1$ at $0$ as a function of
$(\tilde u,\delta,\tilde v)$ and $g_{\alpha}(0,\delta)=0$, 
$\tilde q_{\alpha}(0,\delta,0,h)=0$.
From the triangular structure above, it is obvious that $\tilde{\mathfrak Q}$ is a smooth diffeomorphism at $0\in\R^P$, such that
$\tilde{\mathfrak Q}(0)=0$. Thus, for $\eta<<\epsilon$,
$h\leq h_{0}$ small and $M>>1$, we get the inclusion 
 $\{\vert z_\alpha\vert<\epsilon/M, \vert \tilde v_{\alpha,k}\vert<\eta\}) \subset 
 \tilde{\mathfrak Q}(\{\vert \tilde u_\alpha\vert<\epsilon, \vert \tilde v_{\alpha,k}\vert<\eta\})$.  One has by construction $|\det(D_{(u,v)}\mathfrak Q)\vert=
h^{Q-D}|\det(D_{(\tilde u,\tilde v)}\tilde{\mathfrak Q})\vert$. 
The proof of lemma \ref{lem:diffeo_theta} is complete.

\ep
It is now easy to verify that \eqref{2.20bis} holds true.
One has $\det D_{(u,v)}\hat\kappa_{\xi}=1$ for all $(u,v)\in\R^P$ and for 
$\frac 1 2<\delta<1$, and $0<\eta<<\epsilon <1/2$ 
there exists some numbers  $-1<\alpha_i<\beta_i<1$, $i=1,\ldots, P-D$ depending only on $\epsilon,\eta,\delta$ and such that $\hat\kappa_{\xi}(\Delta_{\epsilon,\eta})$ is contained in the set
$\widehat\Delta_{\epsilon,\eta}=\{(t,s),\;t\in [-\epsilon h,\epsilon h]^D, s\in\prod_{i=1}^{P-D}[\alpha_i h,\beta_i h]\}$. 
Using again the positivity of $g$ and the change of variable $\hat\kappa$, we obtain,
with a constant $c>$ changing from line to line
\be\ba
&h^{-P}\int_{[-h,h]^P}g(\Phi(t))dt\geq h^{-P}\int_{\widehat\Delta_{\epsilon,\eta}}g(\Phi(t))dt\geq h^{-P}\int_{\hat\kappa_{\xi}(\Delta_{\epsilon,\eta})}g(\Phi(t))dt \\
& \geq ch^{-P}\int_{\Delta_{\epsilon,\eta}}g(\Phi\circ\hat\kappa_{\xi}(u,v))dudv=
ch^{-P}\int_{\Delta_{\epsilon,\eta}}g(\Psi(u,v))dudv
\ea\ee
Thanks to  Lemma \ref{lem:diffeo_theta}, we can use the change of variable $\mathfrak Q_{\epsilon,\eta}$ to get
\be\ba
&h^{-P}\int_{\Delta_{\epsilon,\eta}}g(\Psi(u,v))dudv\geq ch^{D-P-Q}\int_{U_{\epsilon,\eta}}
g(e^{\sum_{\alpha\in\A}z_{\alpha}Z^\alpha}w)dzdv \\
&\geq ch^{-Q}\int_{I_{\epsilon',\eta}}g(e^{\sum_{\alpha\in\A}z_{\alpha}Z^\alpha}w)dz
=c h^{-Q}\int_{z\in I_{\epsilon',h}}g(e^{\tilde\lambda (z)}w)\ dz
\ea\ee
whith $\epsilon'=\epsilon/M$ and $M$ is given by Lemma \ref{lem:diffeo_theta}.
The proof of proposition \ref{prop:minor_iter} is complete.
\section{Dirichlet form}\label{sec4}

Let $\mathcal E_{h}$ be the rescaled Dirichlet form associated to the Markov kernel $T_{h}$
\be\label{3.0}
0\leq \mathcal E_{h}(u)= ({1-T_{h}\over h^2}u\vert u)_{L^2}, \quad u\in L^2(M,d\mu)
\ee
The main result of this section is the following proposition.

\begin{prop}\label{prop4.1}
Under the hypoelliptic hypothesis \eqref{hyp_crochet}, there exists $C,h_{0}>0$ such that 
 the following holds true 
for all $h\in ]0,h_{0}]$: 
for all $u\in L^2(M,d\mu)$ such that
\be\label{3.1} 
\Vert u\Vert^2_{L^2} + \mathcal E_{h}(u) \leq 1
\ee
there exists $v_{h}\in \mathcal H^1(\X)$ and $w_{h}\in L^2$ such that 
\be\label{3.2}
u=v_{h}+w_{h}, \quad
\forall j, \ \Vert X_{j} v_{h}\Vert_{L^2} \leq C,  \quad
 \Vert w_{h}\Vert_{L^2}\leq Ch
\ee

\end{prop}

This proposition is easy to prove when the vector fields $X_{j}$ span
the tangent bundle at each point, by elementary Fourier analysis. 
Under the hypoelliptic hypothesis, the proof is more involved,  and will be done in several steps.\\

{\bf Step 1: Localization and reduction to the nilpotent Lie algebra.}

\bigskip
Let us first verify that for all $\varphi\in C^\infty(M)$, there exists $C_{\varphi}$
independent of $h\in ]0,1]$ such that
\be\label{3.300}
\mathcal E_{h}(\varphi u)\leq C_{\varphi}(\Vert u\Vert^2_{L^2}+\mathcal E_{h}(u))
\ee 
One has $1-T_{h}={1\over p}\sum_{k=1}^p(1-T_{k,h})$ and
$$ 2((1-T_{k,h})u\vert u)=\int_{M}{1\over 2h}\int_{-h}^{h}\vert u(x)-u(e^{tX_{k}}x) \vert ^2 dt \ d\mu(x)$$
Since $\sup_{x\in M}\vert \varphi(x)-\varphi(e^{tX_{k}}x)\vert\leq C\vert t \vert$, this implies
for some constant $C_{\varphi}$ and all $k$
$$((1-T_{k,h})\varphi u\vert \varphi u)\leq 
C_{\varphi}(((1-T_{k,h})u\vert u)+h^2\Vert u\Vert^2_{L^2})  $$
and therefore, \eqref{3.300} holds true. Thus, in the proof of proposition
\ref{prop4.1}, we may assume that $u\in L^2(M,d\mu)$ is supported in a small neighborhood 
of a given point $x_{0}\in M$ where theorem \ref{goodman} applies. More precisely,
with the notations of section \ref{sec1}, we may assume in the coordinate system 
$\Lambda\theta$ centered at $x_{0}\simeq 0$ that $u$ is supported in the closed ball
$B^m_{r}=\{x\in \R^m, \vert x \vert \leq r\}\subset V_{0}$. 
Let $\chi(y)\in C_{0}^\infty(U_{0})$
with support in $B^n_{r'}\subset U_{0}$, such that $\int \chi(y)dy=1$. Set $g(x,y)=\chi(y)u(x)$. One has $g(x,y)=
\chi(y)W_{x_{0}}(u)(x,y)$. By hypothesis, one has

$$ \Vert u\Vert^2_{L^2} + \mathcal E_{h}(u) \leq 1$$
This implies for all $k$

$$ 2((1-T_{k,h})u\vert u)=\int_{M}{1\over 2h}\int_{-h}^{h}\vert u(x)-u(e^{tX_{k}}x) \vert ^2 dt \ d\mu(x) \leq ph^2$$
Thus,  for any compact $K\subset U_{0}$, there exist $C_{K}$ such that for all $k$ and
$h\in ]0,h_{0}]$, one has

\be\label{3.4}
\int_{V_{0}\times K}{1\over 2h}\int_{-h}^{h}\vert u(x)-u(e^{tX_{k}}x) \vert ^2 dt \ dxdy \leq C_{K}h^2
\ee
Here, $h_{0}$ is small enough so that $e^{tX_{k}}x$ remains in $V_{0}$ for $\vert t \vert\leq h_{0}$ and $x\in B_{r}$. Let $\phi(x,y)=\chi(y)$. One has $\sup_{x,y}\vert \phi(x,y)-\phi(e^{tZ_{k}}(x,y))\vert\leq C\vert t \vert$ and $\Vert g \Vert_{L^2}\leq C$. Thus,
decreasing $h_{0}$, we get from \eqref{3.4}
that there exists a constant $C$ independent of $k$ and $h\in ]0,h_{0}]$ such that
\be\label{3.5}
\int_{V_{0}\times U_{0}}{1\over 2h}\int_{-h}^{h}\vert g(x,y)-g(e^{tZ_{k}}(x,y)) \vert ^2 dt \ dxdy \leq Ch^2
\ee
Therefore, there exists $C_{0}$ independent  $h\in ]0,h_{0}]$ such that one has 
 \be\label{3.6}
\Vert g \Vert^2_{L^2(\NN)}+ \sum_{j=1}^p h^{-2}
\int_{V_{0}\times U_{0}}{1\over 2h}\int_{-h}^{h}\vert g(x,y)-g(e^{tZ_{k}}(x,y)) \vert ^2 dt \ dxdy
 \leq C_{0}
\ee

\begin{lem}\label{lem4.0}
There exists $C_{1},h_{0}>0$ 
such that for all $h\in ]0,h_{0}]$, any $g$ with support in $B^m_{r}\times B^n_{r'}$,  such that \eqref{3.6} holds true 
can be written on the form 
$$g=f_{h}+l_{h}, \quad \sum_{k=1}^p\Vert Z_{k}f_{h}\Vert_{L^2(V_{0}\times U_{0})}\leq C_{1}, \quad
\Vert l_{h}\Vert_{L^2(V_{0}\times U_{0})}\leq C_{1}h$$
\end{lem}
Let us assume that  lemma \ref{lem4.0} holds true. Then one can write $g=\chi(y)u(x)=f_{h}+l_{h}$.
Let $\psi (x,y)\in C_{0}^\infty(V_{0}\times U_{0})$ equal to $1$ near $B^m_{r}\times B^n_{r'}$.
Set

$$v_{h}=\int \psi(x,y)f_{h}(x,y)dy, \quad w_{h}=\int \psi(x,y)l_{h}(x,y)dy$$
One has $v_{h}+w_{h}=\int \psi(x,y)\chi(y)u(x)dy=\int \chi(y)u(x)dy=u(x)$ and
$\Vert w_{h}\Vert_{L^2}\leq C h$. Moreover, we get from \eqref{2.0}

$$ X_{k}(v_{h})= \int (Z_{k}-\sum_{l}b_{k,l}(x,y){\partial\over \partial y_{l}})\psi(x,y)f_{h}(x,y)dy$$
Since $f_{h}, Z_{k}(f_{h})\in O_{L^2}(1)$ and $\int b{\partial\over \partial y_{l}}(\psi f_{h})dy
=-\int {\partial\over \partial y_{l}}(b)\psi f_{h}dy\in O_{L^2}(1)$, we get that 
\eqref{3.2} holds true. We are thus reduced to prove lemma \ref{lem4.0}.

For any given $k$, the vector field $Z_{k}$ is not singular; thus, decreasing $V_{0}, U_{0}$
if necessary, there exists coordinates $(z_{1},..,z_{D})=(z_{1},z')$ such that 
$Z_{k}={\partial\over \partial z_{1}}$. Using Fourier transform in $z_{1}$, we get that 
if $g$ satisfies  \eqref{3.6}, one has

\be\label{3.7}
2\int (1-{\sin h\zeta_{1}\over h\zeta_{1} })
\vert \hat g(\zeta_{1},z') \vert^2
\ d\zeta_{1}dz'= \int{1\over 2h}\int_{-h}^{h}\vert 1-e^{it\zeta_{1}}\vert ^2 dt 
\vert \hat g(\zeta_{1},z') \vert^2
\ d\zeta_{1}dz'\leq C_{0}'h^2
\ee
Let $a>0$ small. There exists $c>0$ such that 
$(1-{\sin h\zeta_{1}\over h\zeta_{1} })\geq ch^2\zeta_{1}^2$ for $h\vert \zeta_{1}\vert\leq a$
and  $(1-{\sin h\zeta_{1}\over h\zeta_{1} })\geq c$ for $h\vert \zeta_{1}\vert > a$.
Since $g(z_{1},z')={1\over 2\pi}\int_{h\vert \zeta_{1}\vert\leq a}e^{iz_{1}\zeta_{1}}\hat g(\zeta_{1},z')
d\zeta_{1}+{1\over 2\pi}\int_{h\vert\zeta_{1}\vert>a}e^{iz_{1}\zeta_{1}}\hat g(\zeta_{1},z')
d\zeta_{1}= v_{h,k}+w_{h,k}$, we get from \eqref{3.7} that $g$ satisfies
for some $C_{0}$ independent of $h\in ]0,h_{0}]$

\be\label{3.8}
\ba
& \Vert g\Vert_{L^2(\NN)}\leq C_{0}, \quad \text{support}(g)\subset V_{0}\times U_{0}\\
& \forall k, \ \ g=v_{h,k}+w_{h,k}\\
& \Vert Z_{k}v_{h,k}\Vert_{L^2(\NN)}\leq C_{0}, \quad \Vert w_{h,k}\Vert_{L^2(\NN)}\leq C_{0}h
\ea
\ee
and we want to prove that the decomposition $g=v_{h,k}+w_{h,k}$
may be choosen independent of $k$, i.e there exists $C>0$ independent of $h$ such that

\be\label{3.8bis}
\ba
& g=v_{h}+w_{h}\\
& \forall k, \ \Vert Z_{k} v_{h}\Vert_{L^2(\NN)} \leq C \\
& \Vert w_{h}\Vert_{L^2(\NN)}\leq Ch
\ea\ee

In order to prove the implication (\ref{3.8}) $\Rightarrow$ (\ref{3.8bis})
we will construct operators $\Phi, C_{j}, B_{k,j}, R_{l}$, depending on $h$, 
acting on $L^2$ functions with support in a small neighborhood of $o_{\NN}$ in $\NN$, 
with values in $L^2(\NN)$,
such that $\Phi, C_{j}, B_{k,j}, R_{l}, C_{j}hZ_{j}, B_{k,j}hZ_{k} $ 
are uniformly in $h$ bounded on $ L^2$ and 
\be\label{3.9}
\ba
& 1-\Phi = \sum_{j=1}^p C_{j}hZ_{j} +hR_{0}\\
&Z_{j}\Phi = \sum_{k=1}^p B_{k,j}Z_{k} +R_{j}
\ea\ee
and then we set

$$ v_{h}= \Phi (g), \ \  w_{h}=(1-\Phi)(g)$$
With this decomposition of $g$, we get
$$w_{h}=\sum_{j=1}^p C_{j}hZ_{j}(v_{h,j}+w_{h,j}) +hR_{0}(g)\in O_{L^2}(h)$$
$$Z_{k}(v_{h})= \sum_{j=1}^p B_{j,k}Z_{j}(v_{h,j}+h{1\over h}w_{h,j}) +R_{k}(g)\in O_{L^2}(1)$$

We are thus reduced to prove the existence of the operators $\Phi, C_{j}, B_{k,j},R_{l}$,
with the suitable bounds on $L^2$, and such that \eqref{3.9} holds true. This is a problem on the Lie algebra $\NN$ with vector fields $Z_{j}$ given by the Rothschild-Stein Goodman theorem
\ref{goodman}. We will first do this construction in the special case where
the vector fields $Z_{j}$ are equal to the left invariant vector fields $\tilde Y_{j}$
on $\NN$. In that special case, we will have $R_{l}=0$ in formula
\eqref{3.9}. We will conclude in the general case by a suitable h-pseudodifferential
calculus.

\bigskip
{\bf Step 2: The case of left invariant vector fields on $\NN$.} \\

Let $f\ast u$ be the convolution on $\NN$
$$ f\ast u (x)= \int_{\NN}f(x.y^{-1})u(y)dy =\int_{\NN}f(z)u(z^{-1}.x)dz$$
Here, $dy$ is the left (and right) invariant Haar measure on $\NN$, which is simply equal to the Lebesgue measure $dy_{1}...dy_{\r}$ in the coordinates used in formula \eqref{3.10}.
Then for $u\in L^1(\NN)$, the map $f\mapsto f\ast u$ is bounded on $L^q(\NN)$ by 
$\Vert u \Vert_{L^1}$ for any $q\in [1,\infty]$. The vector fields $\tilde Y_{j}$ are divergence free for the Haar measure $dy$. \\

If $f$ is a function on $\NN$, and $a\in \NN$, let 
$\tau_{a}(f)$ be the function defined by $\tau_{a}(f)(x)=f(a^{-1}.x)$. One has
for any $a\in \NN$ and $Y\in T_{e}\NN\simeq \NN$, $\tau_{a}\tilde Y=\tilde Y\tau_{a}$ and 
the following formula holds true:
\be\label{3.11}\ba
&\tau_{a}(f)= \delta_{a}\ast f\\
& \tilde Yf= f\ast \tilde Y\delta_{e}
\ea\ee

Let us denote by $\mathcal T_{h}$ the scaling  operator $\mathcal T_{h}(f)(x)=h^{-Q}f(h^{-1}.x)$. One has $h.(x^{-1})=(h.x)^{-1}$ and $\mathcal T_{h}(f\ast g)=\mathcal T_{h}(f)\ast\mathcal T_{h}(g)$. The action of $\mathcal T_{h}$ on the space $\mathcal D'(\NN)$ 
of distributions on $\NN$, compatible with the action on functions, is given by
$<\mathcal T_{h}(T),\phi>=<T, x\mapsto \phi(h.x)>$. Thus one has 
$\mathcal T_{h}\delta_{e}=\delta_{e}$ and $\mathcal T_{h}(\tilde Y_{j}(\delta_{e}))=
h\tilde Y_{j}(\delta_{e})$ for $j\in \{1,...,p\}$. \\

Let $\mathcal S (\NN)$ be the Schwartz space on $\NN$, and $\varphi\in \mathcal S (\NN)$, with $\int_{\NN}\varphi(x)dx=1$.
For $h\in ]0,1]$, let $\Phi_{h}$ be the operator defined by
\be
\Phi_{h}(f)=f\ast \varphi_{h}, \ \ \varphi_{h}(x)=h^{-Q}\varphi (h^{-1}.x)=\mathcal T_{h}(\varphi)
\ee
Since the Jacobian of the transformation $x\mapsto h.x$ is equal to $h^Q$,
one has $\Vert \varphi_{h}\Vert_{L^1}=\Vert \varphi \Vert_{L^1}$ for all $h\in ]0,1]$,
and therefore the operator $\Phi_{h}$ is uniformly bounded on $L^2$.\\

If we define the operators $B_{k,j,h}$ by $B_{k,j,h}(f)=f\ast \mathcal T_{h}
(\varphi_{k,j})$,
with $\varphi_{k,j}\in \mathcal S (\NN)$, 
the equation
$$ \tilde Y_{j}\Phi_{h}= \sum_{k=1}^p B_{k,j,h}\tilde Y_{k}$$
is equivalent to find the $\varphi_{k,j} \in \mathcal S(\NN)$ such that
\be\label{3.12}
\tilde Y_{j}\varphi=\sum_{k=1}^p\tilde Y_{k}\delta_{e} \ast \varphi_{k,j}
\ee
One has $\int_{\NN}\tilde Y_{j}(\varphi)(x)dx=0$, and since $f\mapsto \tilde Y_{k}\delta_{e} \ast f$ is the right invariant vector field $\mathcal Z_{k}$ on 
$\NN$ such that $\mathcal Z_{k}(o_{\NN})=Y_{k}$, the equation \eqref{3.12} is solvable thanks to
lemma \ref{lemcoho} of the appendix. Moreover, the operators $\Phi_{h}$, $B_{k,j,h}$
and $B_{k,j,h}h\tilde Y_{k}$ are uniformly in $h\in ]0,1]$ bounded on $L^2$.
(one has 
$B_{k,j,h}(h\tilde Y_{k}(f))=f\ast \mathcal T_{h}(\tilde Y_{k}(\delta_{e})\ast \varphi_{k,j})$
and $\tilde Y_{k}(\delta_{e})\ast \varphi_{k,j}\in \mathcal S(\NN)$).\\
 
Let now $c_{j}\in C^\infty(\NN \setminus \{o_{\NN}\})$ be Schwartz for $\bg x \bg \geq 1$,
and  quasi homogeneous of degree $-Q+1$ near $o_{\NN}$ (i.e $c_{j}(t.x)=t^{-Q+1}c_{j}(x)$
for $0<\bg x \bg \leq 1$ and $t>0$ small). Let $C_{j,h}$ be the operators defined by
$C_{j,h}(f)=f\ast \mathcal T_{h}(c_{j})$. Then the equation $1-\Phi_{h}=\sum_{j}C_{j,h}h\tilde Y_{j}$  is
equivalent to solve
\be\label{3.13}
\delta_{e}-\varphi=\sum_{j}\tilde Y_{j}\delta_{e} \ast c_{j}
\ee
In order to solve \eqref{3.13}, we denote by $E\in C^\infty(\NN \setminus \{o_{\NN}\})$
the (unique) fundamental solution,  quasi homogeneous of degree $-Q+2$ on $\NN$ of the 
hypoelliptic equation (for the existence of $E$, we refer to \cite{Fol75}, theorem (2.1), p.172)
$$\delta_{e}=\sum_{j=1}^p \mathcal Z_{j}^2(E), \quad \mathcal Z_{j}(f)=\tilde Y_{j}\delta_{e} \ast f$$ 
Let $\psi\in C_{0}^\infty(\NN)$ with $\psi(x)=1$ near $e=o_{\NN}$. We will choose
$c_{j}$ of the form
\be\label{3.14}
c_{j}=\psi \mathcal Z_{j}(E)-d_{j}, \quad d_{j}\in \mathcal S(\NN)
\ee
Then the equation \eqref{3.13} is equivalent to solve
\be\label{3.15}
\varphi+\sum_{j=1}^p [\mathcal Z_{j},\psi]\mathcal Z_{j}(E)=\varphi_{0}=\sum_{j=1}^p \mathcal Z_{j}(d_{j})
\ee
One has $\varphi_{0}\in \mathcal S(\NN)$ and $\int_{\NN}\varphi_{0}(x)dx=0$
since $\int_{\NN}\varphi (x)dx=1$ and $\int_{\NN}\sum_{j=1}^p [\mathcal Z_{j},\psi]\mathcal Z_{j}(E)dx
=-\int_{\NN}\sum_{j=1}^p \psi \mathcal Z_{j}^2(E)dx=-1$.
Thus, the equation \eqref{3.12} is solvable thanks to
lemma \ref{lemcoho}. Moreover, since $c_{j}\in L^1(\NN)$, the operators
$C_{j,h}$ are uniformly in $h$ bounded on $L^2$. It remains to verify that the operators
$C_{j,h}h\tilde Y_{j}$ are uniformly in $h$ bounded on $L^2$. One has
$C_{j,h}h\tilde Y_{j}(f)= f\ast \mathcal T_{h}(\mathcal Z_{j}(c_{j}))$. Since
$\Vert T_{h}(f)\Vert_{L^2}=h^{-Q/2}\Vert f \Vert_{L^2} $ it is equivalent to
prove that the operator $g\mapsto g\ast \mathcal Z_{j}(c_{j})$ is bounded on $L^2$.
By construction one has $\mathcal Z_{j}(c_{j})=\psi \mathcal Z_{j}^2(E)+l_{j}, l_{j}\in \mathcal S(\NN)$.
With the terminology of \cite{Fol75}, the distribution $Z_{j}^2(E)$
is  homogeneous of degree $0$ (i.e quasi homogeneous of degree $-Q$), thus of the form $\mathcal Z_{j}^2(E)=a_{j}\delta_{e}+f_{j}$
where $f_{j}\in C^\infty(\NN \setminus \{o_{\NN}\})$,  quasi homogeneous of degree
$-Q$ and such that $\int_{b<\vert u\vert<b'}f_{j}(u)du=0$. Thus by 
\cite{Fol75}, proposition 1.9, p.167, the operator $g\mapsto g\ast \mathcal Z_{j}(c_{j})$
is bounded on $L^2$. \\

\bigskip
{\bf Step 3: A suitable $h$-pseudodifferential calculus on $\NN$.} \\

Let $Z^\alpha$ be the smooth vector fields defined in a neighborhood $\Omega$
of $o_{\NN}$ in $\NN$ given by the Goodman theorem \ref{goodman}. 
In this last step, we will finally construct the operators such that \eqref{3.9} holds true.
We first recall the construction of the map $\Theta(a,b)$ which play a crucial role 
in the construction of a parametrix for hypoelliptic operators in \cite{RoSt76}. 
Let us recall that 
$(Y^\alpha_{a}=H_{\alpha}(Y_{1},...,Y_{p})\in T_{e}\NN,\alpha\in \mathcal A)$ is a basis of $T_{e}\NN$.  
For $a\in \NN$ close to $e$ and $u=\sum_{\alpha\in \mathcal A}u_{\alpha}Y^\alpha\in T_{e}\NN$ close to $0$, let $\Lambda(u)=\sum_{\alpha\in \mathcal A}u_{\alpha}Z^\alpha$ and 
\be\label{3.21}
\Phi(a,u)=e^{\Lambda(u)}a
\ee
Clearly, $(a,u)\mapsto (a,\Phi(a,u))$ is a diffeomorphism of a neighborhood
of $(e,0)$ in $\NN\times T_{e}\NN$ onto a neighborhood of $(e,e)$ in $\NN\times\NN$,
and $\Phi(a,0)=a$. We denote by $\Theta (a,b)$ the map defined in a neighborhood
of $(e,e)$ in $\NN\times \NN$ into a neighborhood of $o_{\NN}$ in $\NN\simeq T_{e}\NN$  by
\be\label{3.22}
\Phi(a,\Theta(a,b))=b
\ee
For $b=\Phi(a,u)$, one has 
$\Phi(b,-u)=e^{\Lambda(-u)}(e^{\Lambda(u)}a)=e^{-\Lambda(u)}(e^{\Lambda(u)}a)=a$
. Thus one has the symmetry relation
\be\label{3.20}
\Theta(a,b)=-\Theta(b,a)=\Theta(b,a)^{-1}
\ee

Observe that in the special case $Z_{j}=\tilde Y_{j}$, $\Lambda(u)$ is equal to
the left invariant vector field on $\NN$ such that $\Lambda(u)(o_{\NN})=u$, i.e 
 $\Lambda (u)=\tilde u$  and 
$ \Phi(a,u)=e^{\tilde u}a=a.u$, and this implies in that case
\be\label{3.22bis}
\Theta(a,b)=a^{-1}.b
\ee

Let $\varphi\in \mathcal S (\NN)$, with $\int_{\NN}\varphi(x)dx=1$.
By step 2, there exists functions $\varphi_{k,j}\in \mathcal S (\NN)$,
and $c_{j}\in C^\infty(\NN \setminus \{o_{\NN}\})$, Schwartz for $\bg x \bg \geq 1$,
quasi homogeneous of degree $-Q+1$ near $o_{\NN}$,  
such that the following holds true.

\be\label{3.24} 
\ba
&\tilde Y_{j}(\varphi)=\sum_{k=1}^p \mathcal Z_{k} (\varphi_{k,j})\\
&\delta_{e}-\varphi=\sum_{j}\mathcal Z_{j}(c_{j})
\ea
\ee

 Let $\omega_{0}\subset\subset \omega_{1}$
be small neighborhoods of $o_{\NN}$ such that $\Theta(y,x)$ is well defined for
$(y,x)\in \omega_{0}\times \omega_{1}$, and  $\chi\in C_{0}^\infty(\omega_{1})$
be equal to $1$ in a neighborhood of $\overline \omega_{0}$. 
We define the operators $\Phi_{h}$,
 $B_{k,j,h}$
and $C_{j,h}$ for $1\leq j,k\leq p$ by the formulas

\be\label{3.25} 
\ba  
& \Phi_{h}(f)(x)=\chi(x)\  h^{-Q}\int_{\NN}\varphi(h^{-1}.\Theta(y,x))f(y)dy \\
& B_{k,j,h}(f)(x)=\chi(x)\ h^{-Q}\int_{\NN}\varphi_{k,j}(h^{-1}.\Theta(y,x))f(y)dy\\
& C_{j,h}(f)(x)=\chi(x)\  h^{-Q}\int_{\NN}c_{j}(h^{-1}.\Theta(y,x))f(y)dy
\ea
\ee

All these operators are of the form 
\be\label{3.26}
A_{h}(f)(x)= h^{-Q}\int_{\NN} g(x,h^{-1}.\Theta(y,x))f(y)dy 
\ee
where the function  $g(x,.)$ is smooth in $x$, with compact support $\omega_{1}$, and takes  values in $L^1(\NN)$, i.e 
$\sup_{x\in\omega_{1}}\Vert \partial_{x}^\beta g(x,.)\Vert_{L^1(\NN)}<\infty$ for all $\beta$.
The function $A_{h}(f)$ is well defined for $f\in L^\infty (\NN)$
such that  support$(f)\subset\omega_{0}$ . We have introduce the cutoff $\chi(x)$ just to
have $A_{h}(f)(x)$ defined for all $x\in \NN$, and one has $A_{h}(f)(x)=0$
for all $x\notin \omega_{1}$. 

\begin{lem}\label{lem4.2}
Let $g(x,.)$ be smooth in $x$ with compact support in $\omega_{1}$,  with values in $L^1(\NN)$. 
Then the operator $A_{h}$ defined by \eqref{3.26} is uniformly
in $h\in ]0,1]$ bounded from $L^q(\omega_{0})$ into $L^q(\NN)$ 
for all $q\in [1,\infty]$.
\end{lem}
\bp
The proof is standard. By interpolation, it is sufficient to treat the two cases $q=\infty$ and $q=1$. In the case  $q=\infty$, the jacobian of the change of coordinates 
$y\mapsto u=\Theta (y,x)$ is bounded by $C$ for all $x\in \omega_{1}, y\in \omega_{0}$. Thus we get
$$ \vert A_{h}(f)(x) \vert \leq C \Vert f\Vert_{L^\infty(\omega_{0})}
h^{-Q}\int_{\NN} \vert g(x,h^{-1}.u)\vert du = 
C \Vert f\Vert_{L^\infty(\omega_{0})}\Vert g(x,.)\Vert_{L^1} $$
Since $x\mapsto g(x,.)$ is smooth in $x$ with values in $L^1(\NN)$,
one has $C_{\infty}=\sup_{x\in\omega_{1}}\Vert g(x,.)\Vert_{L^1}<\infty$.
Thus  we get  $\Vert A_{h}(f) \Vert_{L^\infty} \leq CC_{\infty}
\Vert f\Vert_{L^\infty(\omega_{0})}$.\\
For $q=1$, we first extend $g$ as a smooth $L$-periodic function of $x\in\NN$, 
with $L$ large enough,
   $g(x,u)=\sum_{k\in \mathbb Z^D}g_{k}(u)e^{2i\pi k.x/L}$, the equality being  valid for $x\in \omega_{1}$. Observe that  $\Vert g_{k}\Vert_{L^1(\NN)}$ is
   rapidly decreasing in $k$. Then one has 
$$A_{h}(f)(x)=\sum_{k}A_{h,k}(f)(x)e^{ik.x/L}, \quad A_{h,k}(f)(x)=\ h^{-Q}\int_{\NN}
 g_{k}(h^{-1}.\Theta(y,x))f(y)dy $$
The jacobian of the change of coordinates 
$(x,y)\mapsto (u=\Theta (y,x),y)$ is bounded by $C$ for all $(x,y)\in \omega_{1}\times \omega_{0}$, and one has
$$ \int_{\omega_{1}}\vert A_{h,k}(f)(x) \vert dx  \leq C 
h^{-Q}\int_{\NN} \int_{\omega_{0}}\vert g_{k}(h^{-1}.u)\vert \vert f(y)\vert dy du = 
C \Vert f\Vert_{L^1}\Vert g_{k}\Vert_{L^1} $$
Thus we get $\sup_{h\in ]0,1]}\Vert A_{h,k}\Vert_{L^1}= d_{k}$ with $d_{k}$
 rapidly decreasing in $k$, and this implies\\  
 $\sup_{h\in ]0,1]}\Vert A_{h}\Vert_{L^1}\leq \sum_{k}d_{k}<\infty$.
 The proof of lemma \ref{lem4.2} is complete.
\ep

Observe that in the special case $Z_{j}=\tilde Y_{j}$, using \eqref{3.22bis}, we  get that the operators $\Phi_{h}, B_{k,j,h}, C_{j,h}$
defined by the formula \eqref{3.25} are precisely equal, up to the factor $\chi(x)$,
to the operators we have constructed in step 2. \\

In the general case, it remains to show the following:\\

\noindent
i) The operators
$R_{l,h}$ defined by

\be\label{3.9bis}
\ba
& R_{0,h}=h^{-1}\Big (1-\Phi_{h}-\sum_{j=1}^p C_{j,h}hZ_{j}\Big)\\
&R_{j,h}=  Z_{j}\Phi_{h} -\sum_{k=1}^p B_{k,j,h}Z_{k}, \quad 1\leq j\leq p
\ea\ee
are uniformly bounded in $h\in ]0,1]$ on $L^2$.\\

\noindent
ii) The operators $C_{j,h}hZ_{j}$ and $B_{k,j,h}hZ_{k}, k>0$ are
uniformly bounded in $h\in ]0,1]$ on $L^2$.\\

For the verification of i) and ii), we just follow the natural strategy
which is developed in  \cite{RoSt76}.
If $f$ is a function defined near $a\in \NN$, let $\Phi_{a}(f)$
be the function defined near $0$ in $\NN \simeq T_{e}\NN$ by
$ \Phi_{a}(f)(u)=f(\Phi(a,u))$.
The following fundamental lemma is proven  in \cite{RoSt76} (theorem 5) and also in \cite{Go78}
(section 5, "Estimation of the error").

\begin{lem}\label{lem4.1}
For all $j\in \{1,...,p\}$, and  $a\in \NN$ near $e$,    the  vector field
$V_{j,a}$ defined near $0$ in $\NN$ 
\be\label{3.23}
V_{j,a}(g)=\Phi_{a}( Z_{j}(\Phi_{a}^{-1}g))-\tilde Y_{j}(g)
\ee
is of order $\leq 0$ at $0$. If we introduce the system of coordinates $(u_{\alpha})=(u_{l,k})$
with $l(\alpha)=\vert \alpha\vert$ and $1\leq k \leq a_{l}=dim(\NN_{l})$, we thus have

\be\label{3.23bis}
V_{j,a}=\sum_{l=1}^\r \ \sum_{k=1}^{a_{l}}v_{j,l,k}(a,u){\partial\over \partial u_{l,k}} 
\ee
where the functions $v_{j,l,k}(a,u)$ are smooth and satisfy 
$v_{j,l,k}(a,u)\in O(\bg u\bg^l)$.

\end{lem}

Let us denote by $A_{h}[g]$ an operator of the form \eqref{3.26}.
Recall that $g(x,u)$ is  
smooth in $x$ with compact support in $\omega_{1}$, with values in $L^1(\NN)$.
More precisely, we have two cases to consider: a) $g$ is Schwartz in $u$,
and b) $g$ is smooth in $u$ in $\NN \setminus \{o_{\NN}\})$, 
Schwartz for $\bg u \bg \geq 1$, and 
quasi homogeneous of degree $-Q+1$ near $o_{\NN}$.  We have to compute the kernel
of the operators $Z_{j}A_{h}[g]$ and $A_{h}[g]Z_{j}$.\\

We first compute the kernel of $Z_{j}A_{h}(g)$. For any fixed $y$, perform the change of coordinates $x=\Phi_{y}(u)$ so that $\Theta(y,x)=u$.
Denote $Z_{j}^x$ the vector field $Z_{j}$ acting on the variable $x$. Using lemma 
\ref{lem4.1},
we get

\be\label{3.31}\ba
&Z_{j}(A_{h}[g](f))(x)= h^{-Q}\int_{\NN}Z_{j}^x (g(x,h^{-1}.\Theta(y,x))) f(y) dy=\\
& h^{-Q}\int_{\NN} h^{-1}(\tilde Y_{j}^u g)(x,h^{-1}.\Theta(y,x))f(y)dy\\
&+ h^{-Q}\int_{\NN}(Z_{j}^x g)(x,h^{-1}.\Theta(y,x)) f(y)dy\\
&+ \sum_{l=1}^\r \ \sum_{k=1}^{a_{l}}h^{-Q}\int_{\NN} v_{j,l,k}(y,\Theta(y,x))
h^{-l}{\partial g\over \partial u_{l,k}}(x,h^{-1}.\Theta(y,x)) f(y)dy
\ea\ee

By lemma \ref{lem4.2}, the second term in \eqref{3.31} is  uniformly 
bounded in $h\in ]0,1]$, from $L^2(\omega_{0})$ into $L^2(\NN)$. The same holds true for the third term. To see this point, following the proof of lemma \ref{lem4.2}, first write $v_{j,l,k}(y,u)=
\sum_{n} v_{j,l,k,n}(u)e^{2i\pi n.y/L}$, with $v_{j,l,k,n}(u)$ rapidly decreasing in $n$
and $O(\bg u\bg^l)$ near $u=o_{\NN}$. We are then reduce to show that an operator
of the form

$$ R_{h}(f)= h^{-Q}\int_{\NN} h^{-l}G(\Theta(y,x))
{\partial g\over \partial u_{l,k}}(x,h^{-1}.\Theta(y,x)) f(y)dy$$
with $G(u)$ smooth and $G(u)\in O(\bg u\bg^l)$, is uniformly bounded in $h\in ]0,1]$ 
from $L^2(\omega_{0})$ into $L^2(\NN)$ by a constant which depends linearly
on a finite number of derivatives of $G$. Clearly, there exists such a constant $C$ such that $h^{-l}\vert G(\Theta(y,x))\vert \leq C \bg h^{-1}. \Theta(y,x) \bg^l$. 
Thus  the result follows from  the proof of lemma \ref{lem4.2}, since $\bg u\bg ^l{\partial g\over\partial u_{l,k}}(x,u)$ is $L^1$ in $u$ in both case a) and b) (the vector field
$\bg u\bg ^l{\partial \over\partial u_{l,k}}$ is of order $0$). \\

If we denote by $R_{h}$ any operator uniformly bounded on $L^2$, we have thus proven

\be\label{3.31bis}
Z_{j} A_{h}[g] = h^{-1} A_{h}[\tilde Y_{j}^u g] +R_{h}
\ee
 
\bigskip 
Let us now compute the kernel of $A_{h}[g]Z_{j}$. The basic observation is the following identity (recall $u^{-1}=-u$ and $\mathcal Z_{j}(f)=\tilde Y_{j}(\delta_{e})\ast f$ is the right 
invariant vector field such that $\mathcal Z_{j}(0)=Y_{j}$) 

\be\label{3.32}
-\tilde Y_{j}\Big( f(-u)\Big)=\mathcal Z_{j}(f)(-u)
\ee
Let $l_{j}$ be the smooth function such that $^tZ_{j}=-Z_{j}+l_{j}$. For any given $x$
perform the change of coordinates $y=\Phi_{x}(u)$. By \eqref{3.20}, one has
$\Theta(y,x)=-\Theta(x,y)=-u$. We thus get from lemma \ref{lem4.1} and \eqref{3.32}
the following formula:

\be\label{3.33}\ba
&A_{h}[g](Z_{j}(f))(x)= h^{-Q}\int_{\NN}g(x,h^{-1}.\Theta(y,x))Z_{j}(f)(y)dy\\
&= h^{-Q}\int_{\NN}(-Z_{j}^y+l_{j}(y))(g(x,h^{-1}.\Theta(y,x)))f(y)dy\\
& =h^{-Q}\int_{\NN}h^{-1}(\mathcal Z_{j}^u g)(x,h^{-1}.\Theta(y,x))f(y)dy\\
&+h^{-Q}\int_{\NN} g(x,h^{-1}.\Theta(y,x))l_{j}(y)f(y)dy \\
& + \sum_{l=1}^\r \ \sum_{k=1}^{a_{l}} h^{-Q}\int_{\NN} v_{j,l,k}(x,-\Theta(y,x)) h^{-l}
{\partial g \over \partial u_{l,k}}(x,h^{-1}.\Theta(y,x)) f(y) dy
\ea\ee
As above, this gives the identity, with $R_{h}$ uniformly bounded on $L^2$

\be\label{3.33bis}
A_{h}[g]Z_{j}= h^{-1} A_{h}[\mathcal Z_{j}^u g] +R_{h}
\ee

Observe that formulas \eqref{3.24}, \eqref{3.31bis} and \eqref{3.33bis} imply that
\eqref{3.9bis} holds true. Moreover, from \eqref{3.33bis} and lemma
\ref{lem4.2}, the operators  $B_{k,j,h}hZ_{k}, k>0$ are uniformly bounded in $h\in ]0,1]$ on $L^2$. In order to get from \eqref{3.33bis} the same uniform bounds for the operators $C_{j,h}hZ_{j}$, we just observe that in the case  where $g(x,u)$ is quasi homogeneous
in $u$ of degree $-Q+1$ near $o_{\NN}$, one has $\mathcal Z_{j}^u g(x,u)=
C_{j}(x)\delta_{e}+ f_{j}(x,u)$ with $\int_{b<\vert u\vert<b'}f_{j}(x,u)du=0$ and we conclude
as in the end of step 2  by the proposition 1.9
of \cite{Fol75}.\\
The proof of proposition \ref{prop4.1} is complete.

\section{Proof of theorems \ref{thm1} and \ref{thm2}}\label{TV}
This section is devoted to the proof of theorems \ref{thm1} and \ref{thm2}.
Let $\mathcal B_{h}$ be the bilinear form associated to the rescaled 
Dirichlet form $\mathcal E_{h}$
\be
\mathcal B_{h}(f,g)=({1-T_{h}\over h^2} f \vert g)_{L^2}, \quad f,g\in L^2(M,d\mu)
\ee

\begin{prop}\label{propdir}
Let $f\in \mathcal H^1(\X)$. 
Let $(r_{h},\gamma_{h})\in \mathcal H^1(\X)\times L^2$  such that $r_{h}$ converge weakly
(when $h\rightarrow 0$) in $\mathcal H^1(\X)$ to $r\in \mathcal H^1(\X)$, and 
$\sup_{h}\Vert \gamma_{h}\Vert_{L^2}<\infty$. Then
\be\label{5.2}
 \lim_{h\rightarrow 0} \mathcal B_{h}(f,r_{h}+h\gamma_{h})={1\over 6p}\sum_{k=1}^p (X_{k}f \vert X_{k}r)_{L^2}
\ee
\end{prop}
\bp Write $r_{h}=r+r'_{h}$. The weak limit of $r'_{h}$ in $\mathcal H^1(\X)$ is $0$. Since
$\mathcal B_{h}(f,r_{h})=\mathcal B_{h}(f,r)+\mathcal B_{h}(f,r'_{h})$,
we have to prove the two assertions:

\be\label{5.3bis}
\lim_{h\rightarrow 0}\mathcal B_{h}(f,r)=
{1\over 6p}\sum_{k=1}^p (X_{k}f \vert X_{k}r)_{L^2}, \quad \forall f,r \in \mathcal H^1(\X)
\ee
and under the hypothesis that  the weak limit of $r_{h}$ in $\mathcal H^1(\X)$ is $0$
\be\label{5.3}
\lim_{h\rightarrow 0}({1-T_{k,h}\over h^2} f \vert r_{h}+h\gamma_{h})_{L^2}=0,
\quad \forall k\in \{1,...,p\}
\ee

In order to verify \eqref{5.3}, since $M$ is compact, we may assume that $f$ is supported
in a small neighborhood of a point $x_{0}\in M$ where the Goodman theorem \ref{goodman}
applies. 
With the notations of section \ref{sec1},  we may thus assume in the coordinate system 
$\Lambda\theta$ centered at $x_{0}\simeq 0$ that $f,r_{h},\gamma_{h}$ are supported in the closed ball
$B^m_{r}=\{x\in \R^m, \vert x \vert \leq r\}\subset V_{0}$. Let $\chi(y)\in C_{0}^\infty(U_{0})$ with support in $B^n_{r'}\subset U_{0}$, 
such that $\int\chi(y)dy=1$ and write $d\mu(x)=\rho(x)dx$ with $\rho$ smooth.  For $u,v \in L^2(M)$  supported in $B^m_{r}$,
one has 
$$(u\vert v)_{L^2}=\int_{V_{0}} u(x)\overline v (x)\ d\mu(x)=
\int_{V_{0}\times U_{0}} u(x)\overline {\rho(x)\chi(y)v (x)}\ dxdy$$ 
Set $\tilde f(x,y)=W_{x_{0}}(f)(x,y)=f(x), \tilde  r_{h}(x,y)=\rho(x)\chi(y)r_{h}(x),  \tilde \gamma_{h}(x,y)=\rho(y)\chi(y)\gamma_{h}(x)$. 
 We get from \eqref{act_W_T1} 
\be\label{5.4}
({1-T_{k,h}\over h^2} f \vert r_{h}+h\gamma_{h})_{L^2}=
\int_{V_{0}\times U_{0}}\Big({1-\tilde T_{k,h}\over h^2}\tilde f\Big)
\overline{\tilde r_{h}+h \tilde \gamma_{h}}\ dxdy
\ee
Observe that $\tilde \gamma_{h}$ is bounded in $L^2(V_{0}\times U_{0})$.
Since the injection $\mathcal H^1(\X)\subset L^2(M)$ is compact, $r_{h}$ converge
strongly to $0$ in $L^2$, and therefore
$\tilde r_{h}$ converge
strongly to $0$ in $L^2(V_{0}\times U_{0})$. Moreover, 
$Z_{k}(\tilde r_{h})$ converge weakly to $0$ in $L^2(V_{0}\times U_{0})$.
Finally, since $\tilde T_{k,h}$ increase the support of at most
$\simeq h$, we may replace $\tilde f$ by $F=\theta(y) \tilde f$ with $\theta \in C_{0}^\infty$
equal to $1$ near the support of $\chi$. Then $F$ is compactly supported in 
$V_{0}\times U_{0}$ and satisfies $F\in L^2$ and $Z_{k}F\in L^2$. 
Since the vector field $Z_{k}$ is not singular, decreasing $V_{0}, U_{0}$
if necessary, there exists coordinates $(z_{1},..,z_{D})=(z_{1},z')$ such that 
$Z_{k}={\partial\over \partial z_{1}}$. 
One has $dxdy=q(z)dz$ with $q>0$ smooth. Set 
$q\tilde  r_{h}= R_{h}, q\tilde  \gamma_{h}=Q_{h}$. Using Fourier transform in $z_{1}$, it remains to show
\be\label{5.5}
\ba
&\lim_{h\rightarrow 0}\ \ I_{h}=0, \ \ I_{h}=h^{-2}\int  (1-{\sin(h\xi_{1})\over h\xi_{1}})\hat F(\xi_{1},z')
\overline {\hat R_{h}(\xi_{1},z')} d\xi_{1}dz' \\
& \lim_{h\rightarrow 0}\ \ J_{h}=0,\ \ J_{h}= h^{-1}\int  (1-{\sin(h\xi_{1})\over h\xi_{1}})\hat F(\xi_{1},z')
\overline {\hat Q_{h}(\xi_{1},z')} d\xi_{1}dz'
\ea\ee
 Recall that
$Q_{h}$ is bounded in $L^2$, $R_{h}$ converge strongly to zero in $L^2$,
$\partial_{z_{1}}R_{h}$ converge weakly to zero in $L^2$ and $F,\partial_{z_{1}}F\in L^2$. 
We write the first integral in \eqref{5.5} on the form
$$I_{h}=\int \psi(h\xi_{1})\xi_{1}\hat  F(\xi_{1},z')
\overline {\xi_{1}\hat R_{h}(\xi_{1},z')} d\xi_{1}dz'$$
with $\psi(x)=x^{-2}(1-{\sin(x)\over x})$. One has $\psi\in C^\infty(\R)$ and
$\vert \psi(x)\vert \leq C{1\over 1+x^2}$. Then we write 
$I_{h}=I_{1,h}+I_{2,h}$ with $I_{1,h}$ defined by the integral over $\vert\xi_{1}\vert \leq M$
and $I_{2,h}$ defined by the integral over $\vert\xi_{1}\vert > M$.
Since $\xi_{1}\hat R_{h}(\xi_{1},z')$ is bounded in $L^2$, and $\psi \in L^\infty$ we get by Cauchy-Schwarz
$$ \vert I_{2,h}\vert \leq C (\int_{\vert\xi_{1}\vert > M} \vert \xi_{1} \hat F(\xi_{1},z')\vert^2 d\xi_{1}dz')^{1/2}
\rightarrow 0 \quad \text{when} \quad M\rightarrow \infty$$
On the other hand, one has $\psi(x)=\psi(0)+ \tau(x)$ with $\psi(0)=1/6$
and $\sup_{x\in \R}\tau(x)/x \leq C_{0}$. 
Thus we get
\be \label{5.6}
 I_{1,h}={1\over 6}\int_{\vert\xi_{1}\vert \leq M} \xi_{1} \hat  F(\xi_{1},z')
\overline {\xi_{1}\hat R_{h}(\xi_{1},z')} d\xi_{1}dz'+
\int_{\vert\xi_{1}\vert \leq M} \tau(h\xi_{1})\xi_{1} \hat  F(\xi_{1},z')
\overline {\xi_{1}\hat R_{h}(\xi_{1},z')} d\xi_{1}dz'
\ee
For any fixed $M$, the first term in \eqref{5.6} goes to $0$ when $h\rightarrow 0$
since $\xi_{1}\hat R_{h}(\xi_{1},z')$ converge weakly to $0$ in $L^2$ and  
$\xi_{1} \hat F(\xi_{1},z')\in L^2$. Since $\xi_{1}\hat R_{h}(\xi_{1},z')$ is bounded in 
$L^2$ by say $A$, by Cauchy-Schwarz, the  second term is bounded by 
$C_{0}hMA\Vert \partial_{z_{1}}F\Vert_{L^2}$. Thus one has
$\lim_{h\rightarrow 0}\ I_{h}=0$.\\
We proceed exactly in the same way to prove $\lim_{h\rightarrow 0}\ J_{h}=0$:
one has with $x\psi=\phi$
$$J_{h}=\int \phi(h\xi_{1})\xi_{1} \hat  F(\xi_{1},z')
\overline {\hat Q_{h}(\xi_{1},z')} d\xi_{1}dz'$$
and we use the fact that $\phi\in L^\infty$, $\hat Q_{h}(\xi_{1},z')$
is bounded in $L^2$, $\phi(0)=0$ and $\phi(x)/x\in L^\infty(\R)$ .

Let us now verify \eqref{5.3bis}. From \eqref{idel}
this is obvious if $f$ is smooth and $r\in \mathcal H^1(\X)$. 
Standard smoothing arguments show that 
$C^\infty (M)$ is dense in $\mathcal H^1(\X)$. Let now $f\in \mathcal H^1(\X)$ and 
choose  $f_{h}\in C^\infty (M)$ converging
strongly to $f$ in $\mathcal H^1(\X)$. Then $\lim_{h\rightarrow 0}(X_{k}f_{h} \vert X_{k}r)_{L^2}=(X_{k}f \vert X_{k}r)_{L^2}$ and from \eqref{5.3} one has also
$\lim_{h\rightarrow 0}\mathcal B_{h}(f_{h},r)=\lim_{h\rightarrow 0}\mathcal B_{h}(r,f_{h})=\mathcal B_{h}(f,r)$.\\
The proof of proposition \ref{propdir} is complete.
\ep

\subsection {Proof of theorem \ref{thm1}.}

Let $\vert \triangle_{h}\vert $ be the rescaled (non negative)
Laplacien associated to the Markov kernel $T_{h}$:
\be\label{5.-1}
\vert \triangle_{h} \vert ={1-T_{h}\over h^2}
\ee

 From proposition \ref{prop4.1} and lemma \ref{lemA1}, 
there exists $h_{0}>0$ and $C_{4}, C_{5}>0$ independent of $h\in ]0,h_{0}]$, such that
$Spec(\vert \triangle_{h}\vert)\cap [0,\lambda]$ is discrete for all $\lambda \leq C_{4}h^{-2}$ and 
one has the Weyl type estimate 
\be\label{5.0}
\quad \# ( Spec(\vert \triangle_{h}\vert)\cap [0,\lambda] ) \leq C_{5}<\lambda>^{dim(M)/2s}, \quad
\forall \lambda \leq C_{4}h^{-2}. 
\ee 
In particular, since $T_{h}(1)=1$, $1$ is an isolated eigenvalue of $T_{h}$.
Let us verify that $1$ is a simple eigenvalue of $T_{h}$.
Let $f\in L^2=L^2(M,d\mu)$ such that $T_{h}(f)=f$. One has for any $g\in L^2$
\be\label{5.1}
 ((1-T_{h})g\vert g)_{L^2}={1\over 2}\int\int \vert g(x)-g(y)\vert^2\ t_{h}(x,dy)d\mu(x)
\ee
Thus we get for all $k\in \{1,...,p\}$
$$ \int_{M}\int_{-h}^h \vert f(x)-f(e^{tX_{k}}x)\vert^2\ dtd\mu(x)=0$$
This gives $f(x)-f(e^{tX_{k}}x)=0$ for almost all $(x,h)\in M\times ]-h,h[$.
Therefore, one has $X_{k}f=0$ in $\mathcal D'(M)$ for all $k$, and this implies
$f=Cte$ thanks to H\"ormander and Chow theorems. One can also give a
more direct argument: one has $T_{h}^P(f)=f$, and therefore if one use
\ref{5.1} with the Markov kernel $T_{h}^P$ and proposition \ref{prop:minor_iter}, we get

$$\int_{M}\int_{u\in I_{\epsilon,h}}\vert f(x)-f(e^{\lambda (u)}x)\vert^2\ dud\mu(x)=0$$
Since $u\mapsto e^{\lambda (u)}x$ is a submersion, this implies
$f(x)-f(y)=0$ for almost all $(x,y)$ in a neighborhood of the diagonal in $M\times M$,
and therefore $f=Cte$.

Let us now verify that there exists $\delta_{1}>0$ such that for all 
$h\in ]0,h_{0}]$, the spectrum of $T_{h}$ is a subset of $[-1+\delta_{1}, 1]$.
It is sufficient to prove that the same holds true for an odd power $T_{h}^{2N+1}$ of $T_{h}$.
We are thus reduce to show that there exists $h_0, C_0>0$
such that the following inequality holds true for all $h\in ]0,h_0]$
and all $f\in L^2(\Omega)$:
\begin{equation}\label{inf2}
(f+T_{h}^{2N+1}f\vert f)_{L^2}={1\over 2}
\int_{M\times M}t_{h}^{2N+1}(x,dy)\vert f(x)+f(y)\vert^2 d\mu(x)
\geq C_0\Vert f\Vert ^2_{L^2}.
\end{equation}
Take $N$ large enough such that  proposition \ref{prop:minor_iter} applies for 
$T_{h}^{2N+1}$, i.e $t_{h}^{2N+1}(x,dy)\geq c S^\epsilon_{h}(x,dy)$.
Then we are reduce to show that there exists $C$ independent of $h$ such that
\begin{equation}\label{inf1}
\int_{M\times M}S^\epsilon_{h}(x,dy) \vert f(x)+f(y)\vert^2 d\mu(x)
\geq C\Vert f\Vert ^2_{L^2}.
\end{equation}
From the definition  \eqref{def_Sj} of $S^\epsilon_{h}$, we get
$$\int_{M\times M}S^\epsilon_{h}(x,dy) \vert f(x)+f(y)\vert^2 d\mu(x)=
\int_{M}h^{-Q}\int_{u\in I_{\epsilon,h}}\vert f(x)+f(e^{\lambda(u)}x) \vert^2 dud\mu(x)=B$$
Define $A$ by the formula
$$ A=\int_{M}h^{-2Q}\int_{u\in I_{\epsilon/2,h}}\int_{v\in I_{\epsilon/2,h}}
\vert f(e^{\lambda(v)}y)+f(e^{\lambda(u)}y) \vert^2 dudvd\mu(y)$$
Since $\lambda(v)$ is divergence free as a linear combination with
constant coefficients of commutators of the vector fields $X_{k}$, the change of variables $e^{\lambda(v)}y=x$ gives
$$ A=\int_{M}h^{-2Q}\int_{u\in I_{\epsilon/2,h}}\int_{v\in I_{\epsilon/2,h}}
\vert f(x)+f(e^{\lambda(u-v)}x) \vert^2 dudvd\mu(x)$$
Therefore, one has for some constant $c_{\epsilon}>0$ independent of $h$, 
$B\geq c_{\epsilon}A$.
Clearly, one has
$$\int_{M}Re\Big(\int_{u\in I_{\epsilon/2,h}}\int_{v\in I_{\epsilon/2,h}}
f(e^{\lambda(v)}y)\overline f(e^{\lambda(u)}y)  dudv\Big)d\mu(y) \geq 0$$
and this implies, still using the change of variables $e^{\lambda(v)}y=x$
\be\label{inf3}
\ba
&A\geq  2\int_{M}h^{-2Q}\int_{u\in I_{\epsilon/2,h}}\int_{v\in I_{\epsilon/2,h}}
\vert f(e^{\lambda(v)}y)\vert^2 dudvd\mu(y)\\
&=
2\epsilon^D\int_{M}h^{-Q}\int_{v\in I_{\epsilon/2,h}}
\vert f(e^{\lambda(v)}y)\vert^2 dvd\mu(y)=2\epsilon^{2D}\int_{M}\vert f(x)\vert^2d\mu(x)
\ea\ee
From \eqref{inf3} and $B\geq c_{\epsilon}A$, we get that \eqref{inf1} holds true.

\begin{lem}\label{lem:gap}
There exists $C_{2},C_{3}>0$ such that the spectral gap of $T_{h}$ satisfies
\be\label{3.3}
C_{2}h^2 \leq g(h)\leq C_{3}h^2
\ee
\end{lem}
\bp 
The right inequality in \eqref{3.3} is an obvious consequence of the min-max principle
since for any $f\in C^\infty(M)$ one has
$\lim_{h\rightarrow 0}{1-T_{h}\over h^2}f=L(f)$. From \eqref{5.0}, we get that for any 
$a\in ]0,1]$,
 $m_{a}=\sharp (Spec(T_{h})\cap [1-ah^2,1[)$ is bounded by a constant independent of $h$ small, and we have to verify that there exists
 $h_{0}>0$ and $a>0$ independent of $h\in ]0,h_{0}]$  such that $m_{a}=0$. If this is not true, there exists two sequences $\epsilon_{n},h_{n}\rightarrow 0$ and a sequence $f_{n}\in L^2$, with $\Vert f_{n}\Vert_{L^2}=1$ and $(f_{n}\vert 1)_{L^2}= \int_{M}f_{n}d\mu=0$ such that
 $$T_{h_{n}}f_{n}=(1-h_{n}^2\epsilon_{n})f_{n}$$
 This implies $\mathcal E_{h_{n}}(f_{n})=\epsilon_{n}$. Using
 proposition \ref{prop4.1}, we get $f_{n}=v_{n}+h_{n}\gamma_{n}$ with $\sup_{n}\Vert \gamma_{n}\Vert_{L^2}<\infty$
 and $\Vert v_{n}\Vert_{\mathcal H^1(\X)}\leq C$. The hypoelliptic theorem of H\"ormander implies the existence
of $s>0$ such that one has 
$\mathcal H^1(\X)\subset H^s(M)$, hence  the injection $\mathcal H^1(\X)\subset L^2(M)$
is compact.
As a direct byproduct, we get  (up to extraction of a subsequence) that 
the sequence $f_{n}$ converge strongly in $L^2$ to some $f\in \mathcal H^1(\X)$,
and $v_{n}$ converge weakly in $\mathcal H^1(\X)$ to $f$. Set $v_{n}=f+r_{n}$.
Then $r_{n}$ converge weakly to $0$ in $\mathcal H^1(\X)$, $f_{n}=f+r_{n}+h_{n}\gamma_{n}$, and one has
$$\mathcal E_{h_{n}}(f_{n})=\mathcal E_{h_{n}}(f)+ 
2Re (\mathcal B_{h_{n}}(f,r_{n}+h\gamma_{n}))+\mathcal E_{h_{n}}(r_{n}+h_{n}\gamma_{n})$$
Since one has $\mathcal E_{h}(.)\geq 0$, proposition \ref{propdir} implies 
\be\label{5.2}
{1\over 6p}\sum_{k=1}^p \Vert X_{k}f\Vert^2_{L^2}=\lim_{n\rightarrow \infty}\mathcal E_{h_{n}}(f) \leq
\liminf_{n\rightarrow \infty}\mathcal E_{h_{n}}(f_{n})=0
\ee
and therefore $f=Cte$. But since $f_{n}$ converge strongly in $L^2$ 
to $f$, one has $\Vert f\Vert_{L^2}=1$ and $(f\vert 1)_{L^2}= \int_{M}f d\mu=0$. 
This is a contradiction. The proof of lemma \ref{lem:gap} is complete
\ep

To conclude the proof of theorem \ref{thm1}, it remains to prove the 
total variation estimate \eqref{1.7}. Let $\Pi_0$ be the orthogonal projector in $L^2(M,d\mu)$ onto the space of constant functions

\be\label{T1}
\Pi_0(f)(x)= \int_{M} f d\mu
\ee
Then
\be\label{T2}
2sup_{x\in M}\Vert t^n_{h}(x,dy)-\mu\Vert_{TV}= \Vert T_h^n-\Pi_0\Vert_{L^\infty \rightarrow L^\infty}
\ee
Thus, we have to prove that there exist $C_0,h_0$, such that for any $n$ and any $h\in ]0,h_0]$, one has
\be\label{T3}
 \Vert T_h^n-\Pi_0\Vert_{L^\infty \rightarrow L^\infty} \leq C_0e^{-ng(h)}
\ee
Observe that since  $g(h)\simeq h^2$, and 
$\Vert T_h^n-\Pi_0\Vert_{L^\infty \rightarrow L^\infty}\leq 2$,  in the proof of \eqref{T3}, we may assume 
$n\geq Ch^{-2}$ with $C$ large.
Let $E_{h,L}$ be the (finite dimensional) subspace of 
$L^2(M,d\mu)$ span by the eigenvectors  
$e_{j,h}$ of $\vert\triangle_{h}\vert$, associated with eigenvalues $\lambda_{j,h}\leq C_{4}h^{-2}$, with $C_{4}>0$ small enough.  Here, the subscript $L$ means
"low freqencies". Recall from \eqref{5.0} 
$dim(E_{h,L})\leq Ch^{-dim (M)/2s}$. We will denote by $J_{h}$ the set of indices

\be\label{jh}
J_{h}=\{j, \ \lambda_{j,h}\leq C_{4}h^{-2}\}
\ee

\begin{lem}\label{lem:weyl}
There exist  $p>2$ and  $C$ independent of $h\in ]0,h_{0}]$ such that
for all $u\in E_{h,L}$, the following inequality holds true
\be\label{3.5}
\Vert u\Vert^2_{L^p(M)} \leq C(\mathcal E_{h}(u_{})+\Vert u\Vert^2_{L^2})
\ee
\end{lem}
\bp We denote by $C>0$ a constant independent of $h$, changing from line to line.
Let $u\in E_{h,L}$ such that $\mathcal E_{h}(u)+\Vert u\Vert^2_{L^2}\leq 1$.
From proposition \ref{prop4.1}, one has $u=v_{h}+w_{h}$ with 
$\Vert v_{h}\Vert_{\mathcal H^1(\X)}\leq C$ and 
$\Vert w_{h}\Vert_{L^2}\leq Ch$. From the continuous imbedding
$\mathcal H^1(\X)\subset H^s(M)\subset L^q(M)$ with $s>0, q>2, s=dim(M)(1/2-1/q)$,  we get
$$ \Vert v_{h}\Vert_{L^q}\leq C$$
One has $u=\sum_{\lambda_{j,h}\leq C_{4}h^{-2}}z_{j,h}e_{j,h}$
with $\sum_{\lambda_{j,h}\leq C_{4}h^{-2}}\vert z_{j,h}\vert^2\leq 1$.
From corollary \ref{cor:borne_L2_Linf}, one has for $C_{4}>0$ small enough 
$\Vert e_{j,h}\Vert_{L^\infty}\leq Ch^{-Q/2}$. Therefore by Cauchy-Schwarz we get
\be\label{T4}
\Vert u\Vert_{L^\infty}\leq Ch^{-Q/2}(\sum_{\lambda_{j,h}\leq C_{4}h^{-2}}\vert z_{j,h}\vert^2)^{1/2} (dim (E_{h,L}))^{1/2}\leq Ch^{-Q/2-dim(M)/4s}
\ee
From the proof of proposition \ref{prop4.1} (see lemma \ref{lem4.2}),
one has $\Vert v_{h}\Vert_{L^\infty} \leq C \Vert u\Vert_{L^\infty}$. Thus we get
$\Vert w_{h}\Vert_{L^\infty} \leq \Vert u\Vert_{L^\infty}+\Vert v_{h}\Vert_{L^\infty}\leq C h^{-Q/2-dim(M)/4s}$. Since $\Vert w_{h}\Vert_{L^2}\leq Ch$
we get by interpolation that there exists $q'>2$ such that
$$\Vert w_{h}\Vert_{L^{q'}} \leq C$$
Then \eqref{3.5} holds true with $p=\min (q,q')>2$. The proof of lemma \ref{lem:weyl}
is complete.
\ep

We are now ready to prove \eqref{T3}, essentially following 
the strategy of \cite{DiLeMi11}, but with some simplifications.
We  split $T_h$ in 2 pieces, according to the spectral theory.
We write  $T_h-\Pi_0=T_{h,1}+T_{h,2}$ with
\be\label{T5}
T_{h,1}(x,y)=\sum_{\lambda_{1,h}\leq\lambda_{j,h}\leq C_{4}h^{-2}}
(1-h^2\lambda_{j,h})e_{j,h}(x)e_{j,h}(y)
\ee
 One has
$T_h^n-\Pi_0=T_{h,1}^n+T_{h,2}^n$, and we will get the bound
\eqref{T3} for each of the two terms. We start by very rough bounds.
From $\Vert e_{j,h}\Vert_{L^\infty}\leq Ch^{-Q/2}$,
  $\vert (1-h^2\lambda_{j,h})\vert \leq 1$, we get with $A=Q/2+dim(M)/4s$, as in the proof of 
  \eqref{T4} 
with  $C$
independent of $n\geq 1$ and $h$
\begin{equation}\label{T6}
\Vert T_{h,1}^n\Vert_{L^\infty \to L^\infty}\leq \Vert T_{h,1}^n\Vert_{L^2 \to L^\infty}.
\leq Ch^{-A}
\end{equation}
Since $T_h^n$ is bounded by $1$ on $L^{\infty}$, we get from
$T_h^n-\Pi_0=T_{h,1}^n+T_{h,2}^n$ 
\begin{equation}\label{T7}
\Vert T_{h,2}^n\Vert_{L^\infty \to L^\infty}.
\leq Ch^{-A}
\end{equation}
Let $P$ be the integer defined at the beginning of section \ref{sec2}.
Let $M_{h}$ be the Markov operator $M_{h}=T_{h}^P$. Write $n=kP+r$ with $0\leq r <P$. From  proposition \ref{prop:minor_iter} and corollary \ref{cor:borne_spec_ess} one has 
 $M_h=\rho_h+R_h$ with
\begin{equation}\label{T8}
\begin{aligned}
\Vert \rho_{h}\Vert_{L^\infty \to L^\infty}&\leq \gamma<1,\\
\Vert R_{h}\Vert_{L^2 \to L^\infty}&\leq C_0h^{-Q/2}.
\end{aligned}\end{equation}
From this, we deduce that for any $k=1,2,\dots$, one has
$M_h^k=A_{k,h}+B_{k,h}$, with $A_{1,h}=\rho_h, B_{1,h}=R_h$ and the
recurrence relation $A_{k+1,h}=\rho_hA_{k,h},
B_{k+1,h}=\rho_hB_{k,h}+R_hM_h^k$. Thus one gets, since $M_h^k$ is
bounded by $1$ on $L^{2}$,
\begin{equation}\label{T9}\begin{aligned}
\Vert A_{k,h} \Vert_{L^\infty \to L^\infty}&\leq \gamma^k, \\
\Vert B_{k,h}\Vert_{L^2 \to L^\infty}&\leq C_0h^{-Q/2}(1+\gamma+\dots+\gamma^k)\leq 
C_0h^{-Q/2}\big/(1-\gamma).
\end{aligned}\end{equation}
Let $\theta=1-C_{4}<1$ so that $\Vert T_{h,2}\Vert_{L^2 \to
  L^2}\leq \theta$. Then one has
\begin{equation}\label{T9bis}
\Vert T_{h,2}^n\Vert_{L^\infty \to L^2}\leq \Vert T_{h,2}^n\Vert_{L^2
  \to L^2}\leq \theta^n
\end{equation}
For $m\geq 1$, $k\geq 1$, and $0\leq r<P-1$, one gets, using the fact
that $T_{h}$ is bounded by $1$ on $L^\infty$ and \eqref{T7}, \eqref{T9}, and \eqref{T9bis} 
\begin{equation}\label{T10}
\begin{aligned}
\Vert  T_{h,2}^{kP+r+m}\Vert_{L^\infty \to L^\infty}=&
\Vert  T_{h}^rM_{h}^kT_{h,2}^{m}\Vert_{L^\infty \to L^\infty}\leq \Vert  M_h^kT_{h,2}^{m}\Vert_{L^\infty \to L^\infty} \\
&\leq \Vert  A_{k,h}T_{h,2}^{m}\Vert_{L^\infty \to L^\infty} +\Vert B_{k,h} T_{h,2}^{m}\Vert_{L^\infty \to L^\infty} \\
&\leq Ch^{-A}\gamma^k+C_0h^{-Q/2}\theta^m\big/(1-\gamma).
\end{aligned}\end{equation}
Thus we get, that there exists  $C>0$, $\mu>0$, and a large constant $B>>1$ such that  
\begin{equation}\label{T11}
\Vert T_{h,2}^n\Vert_{L^\infty \to L^\infty}
\leq Ce^{-\mu  n}, \qquad \forall h, \quad \forall n\geq B\log(1/h),
\end{equation}
and thus the contribution of $T_{h,2}^n$ is far smaller than the bound
we have to prove in \eqref{T3}.
It remains to study the contribution of $T_{h,1}^n$.

From  lemma \ref{lem:weyl}, using the interpolation inequality
$\Vert u\Vert_{L^2}^2\leq \Vert u\Vert_{L^p}^{{p\over p-1}}\Vert u\Vert_{L^1}^{{p-2\over p-1}}$,
we deduce the Nash inequality, with $1/d=2-4/p>0$ 
\be\label{T20}
\Vert u\Vert_{L^2}^{2+1/d}\leq C(\mathcal E_{h}(u)+\Vert u\Vert^2_{L^2})
\Vert u\Vert_{L^1}^{1/d}, \quad \forall u\in E_{h,L}
\ee
For $\lambda_{j,h}\leq C_{4}h^{-2}$, one has $h^2\lambda_{j,h}\leq 1$, and thus for any 
$u\in E_{h,L}$, one gets $\mathcal E_{h}(u)\leq
 \Vert u\Vert_{L^2}^{2}-\Vert T_hu\Vert_{L^2}^{2}$, thus we get from \ref{T20}
 \be\label{T21}
\Vert u\Vert_{L^2}^{2+1/d}\leq Ch^{-2}((\Vert u\Vert_{L^2}^{2}-\Vert T_hu\Vert_{L^2}^{2}+h^2\Vert u\Vert^2_{L^2})
\Vert u\Vert_{L^1}^{1/d}, \quad \forall u\in E_{h,L}
\ee
From \eqref{T11} and $T_h^n-\Pi_0=T_{h,1}^n+T_{h,2}^n$, we get that there exists $C_2$ such that 
for all $h$ and all $n\geq B\log(1/h)$ one has 
$\Vert T_{1,h}^n\Vert_{L^\infty \rightarrow L^\infty}
\leq C_2$
and thus since $T_{1,h}$ is self adjoint on $L^2$, 
$\Vert T_{1,h}^n\Vert_{L^1 \rightarrow L^1}
\leq C_{2}$.
Fix $p\simeq B\log(1/h)$. Take $g\in L^2$ such that $\Vert g\Vert_{L^1}\leq 1$ and consider the sequence $c_n, n\geq 0$
\be\label{T24}
c_n= \Vert T_{h,1}^{n+p}g\Vert_{L^2}^2 
\ee
Then, $0\leq c_{n+1}\leq c_n$ and from \ref{T21} and $T_{h,1}^{n+p}g \in E_{h,L} $, we get
\be\label{T25}\ba
&c_n^{1+{1\over 2d}}\leq Ch^{-2} (c_n-c_{n+1}+h^2c_n)\Vert T_{h,1}^{n+p}g \Vert_{L^1}^{1/d}\\
&\leq CC_2^{1/d} h^{-2} (c_n-c_{n+1}+h^2c_n)
\ea\ee
Thus there exist 
$A$ which depends only on 
$C, C_2, d$, such that for all $0\leq n \leq h^{-2}$, one has
$c_n\leq ({Ah^{-2}\over 1+n})^{2d}$ (this is the key point in the argument, for a proof
of this estimate, see \cite{DiSa98}).
Thus for all
$0\leq n\leq h^{-2}$ and with $p\simeq B\log(1/h)$ one has 
\be\label{T26}
\Vert T_{h,1}^{n+p}g\Vert_{L^2}\leq ({Ah^{-2}\over 1+n})^d\Vert g\Vert_{L^1} 
\ee
which implies by duality since $T_{1,h}$ is selfadjoint on $L^2$
\be\label{T27}
\Vert T_{h,1}^{n+p}g\Vert_{L^\infty}\leq ({Ah^{-2}\over 1+n})^d\Vert g\Vert_{L^2} 
\ee
Thus
there exist 
$C_0$ ,
such that for $N\simeq h^{-2}$, one has 
\be\label{T27bis}
\Vert T_{h,1}^{N+p}g\Vert_{L^\infty}\leq C_0\Vert g\Vert_{L^2} 
\ee 
and so we get for any $m\geq 0$ and with $N\simeq h^{-2}$
\be\label{T28}
\Vert T_{h,1}^{N+p+m}g\Vert_{L^\infty}\leq C_0(1-h^2\lambda_{1,h})^m\Vert g\Vert_{L^2} 
\ee 
Thus for $n\geq h^{-2}+N+p$, since $h^2\lambda_{1,h}=g(h)$ and $0\leq (1-r)^m\leq e^{-mr}$ for
 $r\in [0,1]$, we get
\be\label{T29}
\Vert T_{h,1}^n\Vert_{L^\infty \rightarrow L^\infty}
\leq C_0e^{-(n-(N+p))g(h)}=C_0e^{(N+p)g(h)}e^{-ng(h)}
\leq C'_0e^{-ng(h)}
\ee
The proof of theorem \ref{thm1} is complete.

\subsection{Proof of theorem \ref{thm2}}
The proof of Theorem \ref{thm2} is exactly the same that the one given in 
\cite{DiLeMi12}. Let $R>0$ be fixed. If
$\nu_h\in[0,R]$ and $u_h\in L^2(M)$ satisfy $| \triangle_h|
u_h=\nu_h u_h$ and $\| u_h\|_{L^2}=1$, then, thanks to proposition
\ref{prop4.1}, $u_h$ can be decomposed as $u_h=v_h+w_h$ with
$\|w_h\|_{L^2}=O(h)$ and $v_h$ bounded in $\mathcal H^1(\X)$. Hence
(extracting a subsequence if necessary) it may be assumed that
$v_h$ weakly converges in $\mathcal H^1(\X)$ to a limit $v$ and that $\nu_h$
converges to a limit $\nu$. Hence $u_{h}$ converge strongly in $L^2$ to $v$.
It now follows from proposition \ref{propdir} that
for any $f\in C^\infty(M)$,
\be\label{T.40}
\nu (f\vert v)=\lim_{h\rightarrow 0} (f\vert \nu_{h}u_{h})= 
\lim_{h\rightarrow 0} (\vert \triangle_{h}\vert (f)\vert u_{h})=\lim_{h\rightarrow 0} \mathcal B_{h}(f,v_{h}+w_{h})={1\over 6p}\sum_{k=1}^p (X_{k}f \vert X_{k}v)_{L^2}=(f\vert Lv)
\ee

Since $f$ is arbitrary, it follows that $(L-\nu)v =0$ . By the Weyl type estimate \eqref{5.0}
the number of eigenvalues $| \triangle_h|$ in the interval $[0,R]$ is uniformly 
bounded. Moreover, the dimension of an orthonormal basis is preserved by strong limit. So the above argument   proves that
for any $\epsilon>0$ small, there exists $h_\epsilon>0$ such that for $h\in]0,h_\epsilon]$, one has
\begin{equation}
\label{T41}
Spec(\vert\Delta_h\vert)\cap[0,R]\subset\cup_j[\nu_j-\epsilon,\nu_j+\epsilon]
\end{equation}
and
\begin{equation}
\sharp Spec(|\Delta_h|)\cap[\nu_j-\epsilon,\nu_j+\epsilon]\leq m_j
\label{T42}
\end{equation}
The fact that one has equality in  \eqref{T42} for $\epsilon$ small follows exactly like in the proof
of theorem 2 iii) in \cite{DiLeMi12}: this use only proposition \ref{propdir}, 
the min-max principle and a compactness argument. The proof of theorem \ref{thm2} is complete.

\begin{remk}
Observe that the estimate \eqref{3.3} on the spectral gap is a direct consequence
of theorem \ref{thm2}, and moreover observe that in the proof of theorem \ref{thm2}
we only use proposition \ref{propdir} in the special case $f\in C^\infty(M)$,
and that for $f\in C^\infty(M)$, proposition \ref{propdir} is obvious.
However, we think that the fact that  proposition \ref{propdir} holds true for any function
$f\in \mathcal H^1(\X)$ is interesting by itself, and since it is an easy byproduct of proposition
\ref{prop4.1}, we decide to include it in the paper.

\end{remk}

\subsection{Elementary Fourier Analysis}
We conclude this section by collecting some basic results on the Fourier analysis theory
(uniformly with respect to $h$) associated to the spectral decomposition of $T_{h}$. These results are  
consequences of the preceding estimates. 
We start with the following lemma which gives an honest $L^\infty$ estimate
on the eigenfunction $e_{j,h}\in E_{h,L}$. Recall $<x>=(1+x^2)^{1/2}$.
\begin{lem}\label{lem-linfty}
There exists $C$ independent of $h$ such that for any eigenfunction $e_{j,h}\in E_{h,L}$,
$\Vert e_{j,h}\Vert_{L^2}=1$,
associated to the eigenvalue $1-h^2\lambda_{j,h}$
of $T_{h}$  the following inequality holds true
\be\label{T30}
\Vert e_{j,h}\Vert_{L^\infty}\leq C <\lambda_{j,h}>^d
\ee
\end{lem}
\bp
This is a byproduct of the preceding estimate \eqref{T27}. Apply this inequality to $g=e_{j,h}$. This gives
\be\label{T31}
(1-h^2\lambda_{j,h})^{n+p}\Vert e_{j,h}\Vert_{L^\infty}\leq ({Ah^{-2}\over 1+n})^d
\ee
Thus we get with $n\simeq h^{-2}<\lambda_{j,h}>^{-1}$
\be\label{T32}
\Vert e_{j,h}\Vert_{L^\infty}\leq ({Ah^{-2}\over h^{-2}<\lambda_{j,h}>^{-1} })^d
(1-h^2\lambda_{j,h})^{-h^{-2}<\lambda_{j,h}>^{-1}-B\log(1/h)}\leq C <\lambda_{j,h}>^{d}
\ee
The proof of lemma \ref{lem-linfty} is complete.
\ep

Let $h_{0}>0$ be a small given real number.
We will use the following notations. If $X$ is a Banach space, we denote by
$X_{h}$ the space $L^\infty(]0,h_{0}], X)$, i.e the space of functions $h\mapsto x_{h}$ from $h\in ]0,h_{0}]$ into $X$
such that $\sup_{h\in ]0,h_{0}]}\Vert x_{h}\Vert_{X}<\infty$. For $a\geq 0$, the notation 
$x_{h}\in  O_{X}(h^a)$ means that there exists $C$ independent of $h$ such that
$\Vert x_{h}\Vert_{X}\leq Ch^a$, and $x_{h}\in  O_{X}(h^\infty)$
means $x_{h}\in  O_{X}(h^a)$ for all $a$. We denote $C^\infty_{h}=\cap_{k\geq 0}C_{h}^k(M)$.

Let $\Pi_{h,L}$ be the $L^2$-orthogonal projection on $E_{h,L}$, and denote
$\Pi_{h,2}=Id-\Pi_{h,L}$. Let $(e_{j,h})_{j\in J_{h}}$ be an orthonormal basis of $E_{h,L}$
with   $T_{h}(e_{j,h})=(1-h^2\lambda_{j,h})e_{j,h}$. For $f\in L^2$ we denote by 
$c_{j,h}(f)=(f\vert e_{j,h})$ the corresponding
Fourier coefficient of $f$. Recall that $J_{h}$ is defined in \eqref{jh}. 
\begin{prop}\label{propfourier}
Let $f_{h}\in C^\infty_{h}$. For  all integer $N$, the following holds true.\\
\be\label{T50}
\vert\triangle_{h}\vert^Nf_{h}\in C^\infty_{h}\quad \text{and} \quad  
\exists \ C_{N}, \quad  \sup_{h\in ]0,h_{0}]}\sum_{j\in J_{h}}\lambda_{j,h}^N \vert c_{j,h}(f_{h})\vert^2
\leq C_{N}
\ee
Moreover, one has the following estimates
\be\label{T51}
\Pi_{h,L}(f_{h})\in O_{L^\infty(M)}(1)
\ee
\be\label{T52}
\Pi_{h,2}(f_{h})\in O_{L^\infty(M)}(h^N)
\ee
\end{prop}
\bp
Let $X$ be a vector field on $M$, and $f\in C^\infty(M)$. The smooth function 
$F(t,x)=f(e^{tX}x)$ satisfy the transport equation
$$\partial_{t}F=X(f), \quad F(0,x)=f(x)$$
Thus, one has by Taylor expansion at $t=0$, and for any integer $N$
$$ F(t,x)=\sum_{n\leq N}{t^n\over n!}X^n(f)(x)+t^{N+1}r_{N}(t,x)$$
with $r_{N}(t,x)$ smooth.  From the definition of $T_{h}$, we thus get

$$T_{h}f(x)\sum_{n \ \text{even} \ \leq N} {h^n\over (n+1)!}
\Big({1\over p}\sum_{k=1}^pX_{k}^n(f)(x)\Big)+h^{N+1}\tilde r_{N}(h,x)$$
with $\tilde r_{N}(h,x)\in C^\infty_{h}$. This implies for $f_{h}\in C^\infty_{h}$
$$\vert\triangle_{h}\vert f_{h}= L(f_{h}) +h^2g_{h}, \quad g_{h} \in C^\infty_{h}$$
Therefore, one has $\vert\triangle_{h}\vert f_{h}\in C^\infty_{h}$, hence
by induction $\vert\triangle_{h}\vert^N f_{h}\in C^\infty_{h}$ for all $N$. The second assertion of
\eqref{T50} follows from $\sup_{h\in ]0,h_{0}]}\Vert g_{h}\Vert_{L^2}<\infty$
for any $g_{h}\in C^\infty_{h}$ and the fact
$$\sum_{j\in J_{h}}\lambda_{j,h}^N \vert c_{j,h}(f_{h})\vert^2=
\Vert \Pi_{h,L} \vert\triangle_{h}\vert^N f_{h}\Vert^2_{L^2} \leq \Vert  \vert\triangle_{h}\vert^N f_{h}\Vert^2_{L^2}$$
For the proof of \eqref{T51}, we just write
$$\Pi_{h,L}(f_{h})=\sum_{j\in J_{h}}c_{j,h}(f_{h})e_{j,h}$$
and we use the estimate \eqref{T30} of lemma \ref{lem-linfty} to get the bound
$$\Vert \Pi_{h,L}(f_{h})\Vert_{L^\infty}\leq C\sum_{j\in J_{h}}\vert c_{j,h}(f_{h})\vert
<\lambda_{j,h}>^{d}$$
$$\leq C\Big( \sum_{j\in J_{h}}\vert c_{j,h}(f_{h})\vert^2
<\lambda_{j,h}>^{2d+2N}\Big)^{1/2}\Big (\sum_{j\in J_{h}}
<\lambda_{j,h}>^{-2N}\Big)^{1/2}$$
From the Weyl type estimate \eqref{5.0}, there exists $N$ and $C$ independent of $h$
such that 
$$\Big(\sum_{j\in J_{h}}<\lambda_{j,h}>^{-2N}\Big)^{1/2} \leq C$$
and therefore \eqref{T51} follows from \eqref{T50}.
It remains to prove   the estimate \eqref{T52}. We first prove the weaker estimate
\be\label{T52bis}
\Pi_{h,2}(f_{h})\in  O_{L^2(M)}(h^N)
\ee
Observe that $\Pi_{h,2}(f_{h})$ satisfies for all $N\geq 1$ the equation
\be\label{T53}
h^{2N}\Pi_{h,2}(\vert\triangle_{h}\vert^N f_{h})=(h^{2}\vert\triangle_{h}\vert)^N\Pi_{h,2}(f_{h})
=(Id-T_{h}\Pi_{h,2})^N\Pi_{h,2}(f_{h})
\ee
By \eqref{T9bis}, the operator $Id-T_{h}\Pi_{h,2}=Id-T_{h,2}$ is invertible on $L^2$ with inverse bounded by $(1-\theta)^{-1}$. Since $\vert\triangle_{h}\vert^N f_{h} \in C^\infty_{h}$
we get from \eqref{T53} $\Pi_{h,2}(f_{h})\in O_{L^2}(h^{2N})$.\\

 Set $g_{h}=\Pi_{h,2}(f_{h})$.
One has $\vert\triangle_{h}\vert^N f_{h}=\Pi_{h,L}(\vert\triangle_{h}\vert^Nf_{h})+\vert\triangle_{h}\vert^Ng_{h}$. From \eqref{T50} and \eqref{T51}, one has 
$\Pi_{h,L}(\vert\triangle_{h}\vert^Nf_{h})\in O_{L^\infty}(1)$. Thus we get $\vert\triangle_{h}\vert^N g_{h}\in O_{L^\infty}(1)$
for any $N$. Let $M_{h}=T_{h}^P$, 
and $\vert\tilde\triangle_{h}\vert=(Id+ T_{h}+...+T_{h}^{P-1})\vert\triangle_{h}\vert$.
Then $g_{h}$ satisfies the equation 
\be\label{T54}
h^{2}\vert\tilde \triangle_{h}\vert g_{h}
=g_{h}-M_{h}g_{h}
\ee
As in  \eqref{T8}, write $M_{h}=\rho_{h}+ R_{h}$.
Since $T_{h}$ is bounded by $1$ on $L^\infty$, one gets
\be\label{T55}
g_{h}-\rho_{h}g_{h}=h^2r_{h}+ R_{h}g_{h}, \quad 
r_{h}=\vert\tilde \triangle_{h}\vert g_{h} \in O_{L^\infty}(1).
\ee
By the second line of \eqref{T8} and \eqref{T52bis} one has
$R_{h}g_{h}\in O_{L^\infty}(h^\infty)$, and by the first line of \eqref{T8},
the operator $Id-\rho_{h}$ is invertible on $L^\infty$ with inverse bounded 
by $(1-\gamma)^{-1}$. Thus we get from \eqref{T55} $g_{h}\in O_{L^\infty}(h^2)$.
Since $\vert\tilde \triangle_{h}\vert g_{h}=\Pi_{h,2}(\vert\tilde \triangle_{h}\vert f_{h})$
and $\vert\tilde \triangle_{h}\vert f_{h}\in C^\infty_{h}$, the same estimates shows
$\vert\tilde \triangle_{h}\vert g_{h}=r_{h}\in O_{L^\infty}(h^2)$. Then \eqref{T55}
implies $g_{h}\in O_{L^\infty}(h^4)$. By induction, we get $g_{h}\in O_{L^\infty}(h^{2N})$
for all $N$.
The proof of proposition \ref{propfourier} is complete.
\ep

Let $F_{k}=Ker(L-\nu_{k})$. Recall $m_{k}=dim(F_{k})$ is the multiplicity 
of the eigenvalue $\nu_{k}$ of $L$. Let us denote by $\mathcal J_{k}$
the set of indices $j$ such that for $h$ small, $\lambda_{j,h}$ is close to $\nu_{k}$,
and $F_{h,k}=span(e_{j,h}, j\in \mathcal J_{k})$.
By theorem \ref{thm2} and his proof, the set $\mathcal J_{k}$ is independent of
$h\in ]0,h_{k}]$ for $h_{k}$ small, and one has $\sharp(\mathcal J_{k})=dim(F_{h,k})=k$
for $h\in ]0,h_{k}]$ . Let $\Pi_{F_{k}}$ and $\Pi_{F_{h,k}}$ the $L^2$-orthogonal
projectors on $F_{k}$ and $F_{h,k}$.

\begin{lem}\label{lembase}
For all $f\in F_{k}$ one has
\be\label{T56}
\lim_{h\rightarrow 0} \Vert f-\Pi_{F_{h,k}}(f) \Vert_{L^\infty}=0
\ee
\end{lem}
\bp
For $f\in F_{k}$, and $h$ small, one has
\be\label{T57}
f-\Pi_{F_{h,k}}(f)=\sum_{j\in J_{h}\setminus \mathcal J_{k}}c_{j,h}(f)e_{j,h}+\Pi_{h,2}(f)
\ee
One has $f\in C_{h}^\infty$, and thus
by \eqref{T52}, we get 
\be\label{T56ter}
\Pi_{h,2}(f)\in O_{L^\infty}(h^\infty)
\ee
 Since $f\in F_{k}$, for any 
given $j\in J_{h}\setminus \mathcal J_{k}$, one has 
$\lim_{h\rightarrow 0}c_{j,h}(f)=\lim_{h\rightarrow 0}(f\vert e_{j,h})_{L^2}=0$.
Therefore, it remains to proove
\be\label{T56bis}
\lim_{N\rightarrow \infty}\sup_{h\in ]0,h_{0}]}
 \sum_{j\in J_{h}, j\geq N}\vert c_{j,h}(f)\vert \Vert e_{j,h}\Vert_{L^\infty}=0
\ee
Let $N>>\nu_{k}$. From \eqref{T30}, Cauchy-Schwarz inequality, \eqref{T50}, and the Weyl type estimate \eqref{5.0}, there exist $N_{0}$ and a constant $C(f)$ independent of $h$ such that one has the estimate 
\be\label{T58}
\ba
& \sum_{j\in J_{h}, j\geq N}\vert c_{j,h}(f)\vert\Vert e_{j,h}\Vert_{L^\infty}
\leq C\sum_{j\in J_{h}, j\geq N}\vert c_{j,h}(f)\vert
<\lambda_{j,h}>^{d}\\
&\leq C\Big( \sum_{j\in J_{h}}\vert c_{j,h}(f)\vert^2
<\lambda_{j,h}>^{2d+2N_{0}}\Big)^{1/2}\Big (\sum_{j\in J_{h}, j\geq N}
<\lambda_{j,h}>^{-2N_{0}}\Big)^{1/2}\\
&\leq C(f)\sup_{h\in ]0,h_{0}]}\Big (\sum_{j\in J_{h}, j\geq N}
<\lambda_{j,h}>^{-2N_{0}}\Big)^{1/2}\longrightarrow 0\quad (N\rightarrow \infty)
\ea\ee
In fact,  since by \eqref{5.0} one has
$\sharp\{j, \lambda_{j,h}\leq m\}\leq C_{5}<m>^{dim(M)/2s}$, one can choose 
$N_{0}=1+dim(M)/4s$. Then one has
$$\sup_{h\in ]0,h_{0}]}\sum_{j\in J_{h}, j\geq N}
<\lambda_{j,h}>^{-2N_{0}}\leq C_{5}\sum_{m\geq m(N)}<m>^{-2N_{0}}<m+1>^{dim(M)/2s}$$
with $m(N)$ the bigger integer such that $\lambda_{N,h}\geq m(N)$ for any $h\in ]0,h_{0}]$.
Observe that \eqref{5.0} implies $\lim_{N\rightarrow \infty}m(N)=\infty$. 
The proof of lemma \ref{lembase} is complete.
\ep 
\section{The hypoelliptic diffusion}\label{sec6}

We refer to the paper of J.-M. Bismut \cite{B81} and references therein  for a construction of the hypoelliptic diffusion associated to the generator $L$.\\

For a given $x_0\in M$, let $X_{x_0}=\{\omega\in C^0([0,\infty[,M), \ \omega(0)=x_0\}$ be the set of continuous paths from $[0,\infty[$ to $M$, starting at $x_0$,
equipped with the topology of  uniform convergence on 
compact subsets of $[0,\infty[$, and let
$\mathcal B$ be the Borel $\sigma$-field generated by the open sets in $X_{x_0}$. We denote by $W_{x_0}$  the Wiener 
measure on $X_{x_0}$ associated to the hypoelliptic diffusion with generator $L$. Let $p_t(x,y)d\mu(y)$ be the heat kernel, i.e the kernel of the self-adjoint
operator $e^{-tL}, t\geq 0$. Then $W_{x_0}$ is the unique probability on $(X_{x_0},\mathcal B)$, such that for any $0<t_1<t_2<...<t_k$ and any Borel sets $A_1,...,A_k$ in $M$, one has

\be\label{6.0}\ba
&W_{x_0}(\omega(t_1)\in A_1, \omega(t_2)\in A_2,...,\omega(t_k)\in A_k)= \\
&\int_{A_1\times A_2\times ...\times A_k} p_{t_k-t_{k-1}}(x_{k},x_{k-1})...p_{t_2-t_1}(x_2,x_1)p_{t_1}(x_1,x_0)
d\mu(x_1)d\mu(x_2)...d\mu(x_k)
\ea\ee

\bigskip

Let us first introduce some notations. Let $Y=\{1,...,p\}\times [-1,1]$ and let $\rho$
be the uniform probability on $Y$. For any function $g(k,s)$ on $Y$, one has
\be\label{6.1}
\int_{Y}gd\rho={1\over 2p}\sum_{k=1}^p \int_{-1}^{+1}g(k,s)ds
\ee 
We denote by $Y^{\mathbb N}$ the infinite product space 
$ Y^{\mathbb N}= \{\underline y=(y_1,y_2,...,y_n,...), \  y_{j}\in Y\}$.
Equipped with the product topology, it is a compact metrisable space, and we denote by 
$\rho^{\mathbb N}$ the product probability on $Y^{\mathbb N}$.
Let $M^{\mathbb N}$ be the infinit product space
$ M^{\mathbb N}= \{\underline x=(x_1,x_2,...,x_n,...), \  x_{j}\in M\}$.
Equipped with the product topology, $M^{\mathbb N}$ is a compact metrisable space.
For $h\in ]0,1]$, and $x_{0}\in M$, let $\pi_{x_{0},h}$ be the continuous map from 
$Y^{\mathbb N}$ into $M^{\mathbb N}$ defined by
\be\label{6.2}
\pi_{x_{0},h}((k_{j},s_{j})_{j\geq 1})=(x_{j})_{j\geq 1}, \quad 
x_{j}=e^{s_{j}hX_{k_{j}}}...e^{s_{2}hX_{k_{2}}}e^{s_{1}hX_{k_{1}}}x_{0}
\ee
We will use the notation $X^n_{h,x_{0}}=(\pi_{x_{0},h})_{n}$. This means that
$X^n_{h,x_{0}}$ is the position after $n$ step of the random walk starting at $x_{0}$. 
Let $\mathcal P_{x_0,h}$ be the probability on  $M^{\mathbb N}$
defined by $\mathcal P_{x_0,h}=(\pi_{x_{0},h})_{*}(\rho^{\mathbb N})$.
Then by construction, one has for all Borel sets $A_1,...,A_k$ in $M$

\be\label{6.3}
\begin{split}
&\mathcal P_{x_0,h}(x_1\in A_1,x_2\in A_2,...,x_k\in A_k)\\
&=\int_{A_1\times A_2\times ...\times A_k} t_h(x_{k-1},dx_k)...
t_h(x_1,dx_2)t_h(x_0,dx_1)
\end{split}
\ee
Let $j_{x_0,h}$ be the  map
from $Y^{\mathbb N}$ into $X_{x_0}$ defined by, with $\underline y=((k_{j},s_{j})_{j\geq 1})$
\be\label{6.5}\ba
j_{x_0,h}(\underline y)=\omega \iff & \forall j\geq 0, \ \forall t\in [0,h^2], 
\quad \omega(jh^2+t)=e^{{t\over h^2}hs_{j}X_{k_{j}}}x_{j}\\
& \text{with} \ x_{j}=(\pi_{x_{0},h}(\underline y))_{j} \ \text{if} \ j\geq 1
\ea\ee
Let $P_{x_0,h}$ be the probability on $X_{x_0}$ defined as the image
of $\rho^{\mathbb N}$ by the continuous map $j_{x_0,h}$. Our aim is to prove the 
following theorem of weak
convergence of $P_{x_0,h}$
 to the Wiener measure $W_{x_0}$ when $h\rightarrow 0$. 

\begin{thm}\label{th4}
For any bounded continuous function $\omega \mapsto f(\omega)$ on $X_{x_0}$, one has
\be\label{6.6}
\lim_{h\rightarrow 0}\int f dP_{x_0,h}=\int f dW_{x_0}
\ee
\end{thm}

 Observe that the proof below  shows that our study of 
the Markov kernel $T_{h}$
on $M$ is also a way to prove the existence of the
Wiener measure $W_{x_{0}}$ associated to the hypoelliptic diffusion.
Let $g$ be a Riemannian distance on $M$ and let $d_{g}$ the associated distance.
We start by proving that  the family of probability $P_{x_0,h}$ is tight, hence  compact by the Prohorov theorem.\\

\begin{prop}\label{prop6.1}
For any  $\varepsilon>0$, there exists $h_{\varepsilon}>0$ such that the following holds true for any $T>0$. 
\be\label{6.7}
\lim_{\delta\rightarrow 0} \Big(\sup_{h\in ]0,h_{\varepsilon}]}P_{x_0,h}(\max_{\vert s-t\vert\leq\delta, \ 0\leq s,t\leq T}d_g(\omega(s),\omega(t))>\varepsilon )\Big)=0
\ee
\end{prop}
\bp
We start with the following lemma.
\begin{lem}\label{lem6.1}
Let $f\in C^\infty(M)$. There exists $C$  such that
for all $h\in ]0,h_{0}]$, one has
\be\label{6.8}
\forall \delta\in [0,1], \quad \sup_{nh^2\leq\delta}\Vert T_{h}^n(f)-f-nh^2\vert\triangle_{h}\vert f\Vert_{L^\infty}\leq C\delta^2
\ee
\end{lem}
\bp
We may assume $\delta>0$ and $n\geq 1$. Then $nh^2\leq \delta$ implies $h\leq \sqrt\delta$. With the notation of section \ref{TV}, one has
\be\label{6.9}\ba
&T_{h}^n(f)-f-nh^2\vert\triangle_{h}\vert f=
\sum_{j\in J_{h}}c_{j,h}(f)\Big((1-h^2\lambda_{j,h})^n-1-nh^2\lambda_{j,h} \Big)
e_{j,h}+R(n,h)\\
&R(n,h)=T_{h}^n\Pi_{h,2}(f)-\Pi_{h,2}(f+nh^2\vert\triangle_{h}\vert f)
\ea\ee
One has $\vert\triangle_{h}\vert f\in C_{h}^\infty$ by \eqref{T50}, $T_{h}$
is bounded by $1$ on $L^\infty$, and $nh^2\leq\delta\leq1$. Thus from \eqref{T52}
we get
\be\label{6.10}
\sup_{nh^2\leq\delta} \Vert R(n,h)\Vert_{L^\infty}\in O(h^\infty)\subset O(\delta^\infty) 
\ee
For all $j\in J_{h}$ one has $h^2\lambda_{j,h}\in [0,1]$ and for all $x\in [0,1]$
$$\vert (1-x)^n-1-nx\vert\leq{n(n-1)\over 2}x^2$$
Therefore we get
\be\label{6.11}
\Vert \sum_{j\in J_{h}}c_{j,h}(f)\Big((1-h^2\lambda_{j,h})^n-1-nh^2\lambda_{j,h} \Big)e_{j,h}
\Vert_{L^\infty}\leq {n^2h^4\over 2}\sum_{j\in J_{h}}\lambda_{j,h}^2\vert c_{j,h}(f)\vert
 \Vert e_{j,h}\Vert_{L^\infty}
\ee
By the Weyl type estimate \eqref{5.0}, \eqref{T30} and \eqref{T50}, 
there exists a constant $C$ such that 
$$\sup_{h\in ]0,h_{0}]}\sum_{j\in J_{h}}\lambda_{j,h}^2\vert c_{j,h}(f)\vert
 \Vert e_{j,h}\Vert_{L^\infty} \leq C $$
 Therefore \eqref{6.8} is consequence of \eqref{6.10} and \eqref{6.11}.
 The proof of lemma \ref{lem6.1} is complete.
\ep

The proof of proposition \ref{prop6.1} is now standard and proceeds as follows.
Let $\varepsilon_{0} >0$ small with respect to the injectivity radius of the
Riemannian manifold $(M,g)$, and let $\varepsilon\in]0,\varepsilon_{0}]$ be fixed. One has 
\be\label{6.12}
\rho^\mathbb N (d_g(X^n_{h,x_0},x_0)>\varepsilon)=
\int_{d_g(y,x_0)>\varepsilon}t_h^n(x_0,dy)=T_h^n(1_{d_g(y,x_0)>\varepsilon})(x_0)
\ee
Let $\varphi(r) \in C^\infty([0,\infty[)$ be a nondecreasing function 
equal to $0$ for $r\leq 3/4$ and equal to $1$
for $r\geq 1$. Set
\be\label{6.13}
\varphi_{x_0,\varepsilon}(x)=\varphi({d_g(x,x_0)\over\varepsilon}) 
\ee
Then $\varphi_{x_0,\varepsilon}$ is a smooth function , and  from
$1_{d_g(y,x_0)>\varepsilon}\leq \varphi_{x_0,\varepsilon}\leq 1$, we get 
since $T_h$ is Markovian,
\be\label{6.14}
0\leq T_h^n(1_{d_g(y,x_0)>\varepsilon})\leq T_h^n(\varphi_{x_0,\varepsilon})
\ee
Since $T_{h}$ moves the support at distance $\leq ch$, one has
$\varphi_{x_0,\varepsilon}(x_{0})+nh^2(\vert\triangle_{h}\vert \varphi_{x_0,\varepsilon})
(x_{0})=0$
for $ch\leq \varepsilon/2$,
From lemma \ref{lem6.1}, we thus get that there exists $h_{\varepsilon}>0$ 
and $C_{\varepsilon}$
such that
\be\label{6.15}
\sup_{h\in]0,h_{\varepsilon}]}\sup_{nh^2\leq \delta} T_h^n(\varphi_{x_0,\varepsilon})(x_{0})\leq C_{\varepsilon}\delta^2
\ee
Since $M$ is compact, it is clear from the proof of lemma \ref{lem6.1}
that we may assume $C_{\varepsilon}$ independent of $x_{0}\in M$. From \eqref{6.12}, \eqref{6.14} and \eqref{6.15} we get
\be\label{6.16}
\sup_{x_{0}\in M}\sup_{h\in]0,h_{\varepsilon}]}\sup_{nh^2\leq \delta}
\rho^\mathbb N (d_g(X^n_{h,x_0},x_0)>\varepsilon) \leq C_{\varepsilon}\delta^2
\ee
Let $T>0$ be given. One has for  $h\in ]0,h_{\varepsilon}]$ the following inequalities.
\be\label{6.17}\ba
&\rho^\mathbb N(\exists j<l\leq h^{-2}T,(l-j)h^2\leq\delta, \ \ d_g(X^j_{h,x_0},X^l_{h,x_0})>4\varepsilon)\\
\leq &{C\over \delta} \sup_{y_0\in M}\rho^\mathbb N(\exists j<l\leq h^{-2}\delta,\ \ d_g(X^j_{h,y_0},X^l_{h,y_0})>4\varepsilon)\\
\leq &{C\over \delta}\sup_{y_0\in M}\rho^\mathbb N(\exists j\leq h^{-2}\delta, \ \ d_g(X^j_{h,y_0},y_0)>2\varepsilon)\\
\leq &{2C\over \delta} \sup_{z_0\in M, nh^2\leq \delta}\rho^\mathbb N (d_g(X^n_{z_0},z_0)>\varepsilon)\\
\text{(by (\ref{6.16}))}&\leq 2CC_{\varepsilon}\delta
\ea
\ee
In fact, for the first inequality in \eqref{6.17}, we just use the fact that the interval $[0,T]$ is a union
of $\simeq C/\delta$ intervals of length $\delta/2$. The second inequality is obvious since the event
$\{\exists j<l\leq h^{-2}\delta,\ \ d_g(X^j_{h,y_0},X^l_{h,y_0})>4\varepsilon\}$ is a subset of 
$\{\exists j\leq h^{-2}\delta, \ \ d_g(X^j_{h,y_0},y_0)>2\varepsilon\}$. For the third one, we use the fact that the event
$A=\{\exists j\leq h^{-2}\delta, \ \ d_g(X^j_{h,y_0},y_0)>2\varepsilon\}$ is contained in $B\cup_{j<k}(C_j\cap D_j)$
with $B=\{d_g(X^k_{h,y_0},y_0)>\varepsilon\}$ ($k$ is the greatest integer $\leq\delta h^{-2}$), $C_j=\{d_g(X^j_{h,y_0},X^k_{h,y_0})>\varepsilon\}$,
$D_j=\{d_g(X^j_{h,y_0},y_0)>2\varepsilon \ \text{and} \ d_g(X^l_{h,y_0},y_0)\leq 2\varepsilon \ \text{for} \ l<j\}$, 
and the fact that $C_j$ and $D_j$ are independent and the $D_{j}$ are disjoints.\\

Since $P_{x_{0},h}=(j_{x_{0},h})_{*}(\rho^\mathbb N)$, \eqref{6.7} follows easily
from \eqref{6.17} and the definition \eqref{6.5} of the map $j_{x_{0},h}$.
The proof of proposition \ref{prop6.1} is complete.

\ep

With the result of proposition \ref{prop6.1}, the proof of theorem \ref{th4} follows 
now the classical proof of weak convergence
of a sequence of random walks in the Euclidian space $\mathbb R^d$ to 
the Brownian motion on $\mathbb R^d$, for which we refer to (\cite{KS}, chapter 2.4). 
We have to prove that any weak limit $P_{x_0}$ of a sequence $P_{x_0,h_k}$, $h_k\rightarrow 0$, is equal to the Wiener measure $W_{x_0}$. We denote by $\omega_{h}(t)$
the map from $Y^\mathbb N$ into $M$ defined by 
$\omega_{h}(t)(\underline y)=j_{x_{0},h}(\underline y)(t)$. By theorem 4.15 of \cite{KS} it is sufficient  to show that for
any $m\geq 1$, any $0<t_{1}<...<t_{m}$, and any  continuous function $f(x_{1},...,x_{m})$
defined on the space $M^m$,  one has

\be\label{5.120}\ba
&\lim_{h\rightarrow 0}\int_{Y^\mathbb N} f(\omega_{h}(t_1),...,\omega_{h}(t_m))d\rho^\mathbb N=\\
&\int f(x_{1},...,x_{m}) p_{t_m-t_{m-1}}(x_{m},x_{m-1})...p_{t_2-t_1}(x_2,x_1)p_{t_1}(x_1,x_0)
d\mu(x_1)d\mu(x_2)...d\mu(x_m)
\ea\ee
 As in \cite{KS}, we may assume $m=2$. For a given $t\geq 0$,
 let $n(t,h)\in \mathbb N$ be the greatest integer such that 
$h^2n(t,h)\leq t$. By \eqref{6.5}), one has for some $c>0$ independent of $h$ 
and $\underline y\in Y^\mathbb N$, 
$d_{g}(\omega_{h}(t), X^{n(t,h)}_{h,x_0})\leq ch$. Since $f$ is uniformly continuous on
$M^m$,   we are reduce to prove
\be\label{5.121}
\ba
&\lim_{h\rightarrow 0}\int f(X^{n(t_1,h)}_{h,x_0}, X^{n(t_2,h)}_{h,x_0})d\rho^{\mathbb N}\\
&=\int f(x_{1},x_{2}) p_{t_2-t_{1}}(x_{2},x_{1})p_{t_1}(x_1,x_0)
d\mu(x_1)d\mu(x_2)
\ea
\ee
From \eqref{6.3}, one has
\be\label{5.122}
\ba
&\int f(X^{n(t_1,h)}_{h,x_0}, X^{n(t_2,h)}_{h,x_0})d\rho^\mathbb N\\
&=\int f(x_{1},x_{2}) t_h^{n(t_2,h)-n(t_1,h)}(x_1,dx_2)t_h^{n(t_1,h)}(x_0,dx_1)
\ea
\ee
By (\ref{5.121}), (\ref{5.122}), we have to show that for any continuous function
$f(x_1,x_2)$ on the product space $M\times M$, one has
\be\label{5.120a}\ba
&\lim_{h\rightarrow 0}\int_{M\times M}f(x_1,x_2) t_h^{n(t_2,h)-n(t_1,h)}(x_1,dx_2)t_h^{n(t_1,h)}(x_0,dx_1)\\
&=\int_{M\times M}f(x_1,x_2)p_{t_2-t_{1}}(x_{2},x_{1})p_{t_1}(x_1,x_0)
d\mu(x_1)d\mu(x_2)
\ea\ee
or equivalently
\be\label{5.120b}
\ba
\lim_{h\rightarrow 0}T_h^{n(t_1,h)}&\Big(T_h^{n(t_2,h)-n(t_1,h)}(f(x_1,.))(x_1)\Big)
(x_0)\\
&=e^{-t_1L}\Big(e^{-(t_2-t_1)L}(f(x_1,.))(x_1)\Big)
(x_0)
\ea
\ee
Since $\Vert T_h^{n(t,h)} \Vert_{L^\infty}\leq 1$ and
$\Vert e^{-tL} \Vert_{L^\infty}\leq 1$, the following "central limit" theorem will conclude the proof of theorem \ref{th4}. \\

\begin{lem}\label{lem6.2} 
For all $f\in C^0(M)$, and all $t>0$, one has
\be\label{5.103} 
\lim_{h\rightarrow 0}\Vert e^{-tL}(f)-T_h^{n(t,h)}(f)\Vert_{L^\infty} =0
\ee
\end{lem}
Since one has $\Vert T_h^{n(t,h)} \Vert_{L^\infty}\leq 1$ and
$\Vert e^{-tL} \Vert_{L^\infty}\leq 1$, 
it is sufficient to prove that (\ref{5.103}) holds true for $f
\in \mathcal D$, with $\mathcal D$  a dense subset of the space
$C^0(M)$, and therefore we may assume that $f\in F_{k}$ is an eigenvector of $L$
associated to the eigenvalue $\nu_{k}$. We set $n=n(t,h)$,
and we use the notation of section \ref{TV}. One has

\be\label{T57bis}\ba
&T_h^{n}(f)=\sum_{j\in \mathcal J_{k}}c_{j,h}(f)(1-h^2\lambda_{j,h})^n e_{j,h}+ R_{t,h}(f)\\
&R_{t,h}(f)=\sum_{j\in J_{h}\setminus \mathcal J_{k}}c_{j,h}(f)(1-h^2\lambda_{j,h})^n e_{j,h}+
T_h^{n}\Pi_{h,2}(f)
\ea\ee
One has  $\vert (1-h^2\lambda_{j,h})^n\vert\leq 1$ and $T_{h}$ is bounded by $1$ on $L^\infty$.
By \eqref{T56ter} and \eqref{T56bis},  we thus get
$$\lim_{h\rightarrow 0}\Vert R_{t,h}(f) \Vert_{L^\infty}=0$$
One has $\lim_{h\rightarrow 0}(1-h^2\lambda_{j,h})^{n(t,h)}=e^{-t\nu_{k}}$
for all $j\in \mathcal J_{k}$. Moreover, one has $\sharp \mathcal J_{k}=m_{k}$
and $\sup_{h\in ]0,h_{0}]}\sup_{j\in \mathcal J_{k}}\Vert e_{j,h}\Vert_{L^\infty}<\infty$
by lemma \ref{lem-linfty}.
Therefore lemma \ref{lembase} and $e^{-tL}(f)=e^{-t\nu_{k}}f$ implies
$$\lim_{h\rightarrow 0}\Vert \sum_{j\in \mathcal J_{k}}c_{j,h}(f)(1-h^2\lambda_{j,h})^n e_{j,h}
-e^{-tL}(f)\Vert_{L^\infty}=0$$
The proof of lemma \ref{lem6.2} is complete.
\ep

\section{Appendix}

Let $P=P(x,\partial_{x})$ be an elliptic second order differential operator on $M$,
with smooth coefficients,
such that $P=P^*\geq Id$, where $P^*$ is the formal adjoint on $L^2(M,\mu)=L^2$. 
Let $(e_{j})_{j\geq 1}$ be  an
orthonormal basis of eigenfunctions of $P$ in $L^2$,  and $1\leq \nu_{1}\leq \nu_{2}...$ be the associated eigenvalues. By the classical Weyl formula, one has
\be\label{a0}
\# \{j, \ \nu_{j}^{1/2}\leq r\} \simeq r^{dim(X)}
\ee
For $s\in \mathbb R$ and $f=\sum_{j}f_{j}e_{j}$ in the Sobolev space $H^s(M)$, we set
$$\Vert v\Vert_{H^s}^2=\sum_{j} \nu_{j}^s\vert f_{j}\vert^2=(P^s f\vert f)_{L^2}$$
Let us recall that this $H^s$-norm depends on $P$, but an  other choice for $P$ 
gives an equivalent norm. The following elementary lemma is useful for us.  

\begin{lem}\label{lemA1}
Let $s>0$ and $A_{h}=A^*_{h}\geq 0$, $h\in ]0,1]$ 
a family of non negative self-adjoint bounded operators
acting on $L^2(M,\mu)$. Assume that there
exists a constant $C_{0}>0$ independent of $h$ such that 
\be \label{a1}
((Id+A_{h})u\vert u) \leq 1 \Rightarrow \exists (v,w)\in H^s\times L^2 \ \text{such that}  \ u=v+w, \ \Vert v\Vert_{H^s}\leq C_{0}, 
\ \Vert w\Vert_{L^2}\leq C_{0}h 
\ee
Let $C_{1}<{1\over 4C_{0}^2}$. There exists $C_{2}>0$ independent of $h$ such that
$Spec(A_{h})\cap [0,\lambda-1]$ is discrete for all $\lambda \leq C_{1}h^{-2}$ and 
\be\label{a2}
\quad \# ( Spec(A_{h})\cap [0,\lambda-1] ) \leq C_{2}<\lambda>^{dim(M)/2s}, \quad
\forall \lambda \leq C_{1}h^{-2}
\ee 
Here, $\# ( Spec(A_{h})\cap [0,r])$ is the number of eigenvalues of $A_{h}$
in the interval $[0,r]$ with multiplicities, and $<\lambda>=\sqrt{1+\lambda^2}$.
\end{lem}
\bp
Let $B_{h}=Id+A_{h}$. Let $C_{h}$ be the bounded operator on $L^2$ defined by
$$ C_{h}(\sum_{j}u_{j}e_{j})=\sum_{j}\min (h^{-1}, \nu_{j}^{s/2})u_{j}e_{j}$$
For $u=v+w$ one has
$$ \Vert C_{h}u\Vert_{L^2}^2\leq 2 \Vert C_{h}v\Vert_{L^2}^2+2\Vert C_{h}w\Vert_{L^2}^2
\leq 2(\Vert v\Vert_{H^s}^2+h^{-2}\Vert w \Vert_{L^2}^2)$$
From \eqref{a1}, we get for all $u\in L^2$
\be\label{a3}
\Vert C_{h}u\Vert_{L^2}^2\leq 4C_{0}^2(B_{h}u\vert u)
\ee
For any non negative selfadjoint bounded operator $T$ on $L^2$, set for $j\geq 1$
$$\lambda_{j}(T)=\min_{dim(F)=j}(\max_{u\in F, \Vert u\Vert_{L^2}=1}(Tu\vert u))$$
It is well known that if $\#\{j, \lambda_{j}(T)\in [0,a[\}<\infty$, the spectrum of
$T$ in $[0,a[$ is discrete and in that case, the $\lambda_{j}(T)\in [0,a[$ are the eigenvalues
of $T$ in $[0,a[$ with multiplicities. 
From \eqref{a3}, we get for all $j\geq 1$ the inequality
\be\label{a4}
\lambda_{j}(B_{h})\geq {1\over 4C_{0}^2}\lambda_{j}(C_{h}^2)
\ee
For all $j$ such that $\nu_{j}^s<h^{-2}$, one has $\lambda_{j}(C_{h}^2)=\nu_{j}^s$,
and therefore, for all $\lambda<h^{-2}$, we get from  \eqref{a0}, 
 $\#\{j, \lambda_{j}(C_{h}^2)\leq \lambda \}\leq C
<\lambda>^{dim(M)/2s}$. Therefore, the spectrum of $B_{h}$ in $[0,h^{-2}/{4C_{0}^2}[$
is discrete, and \eqref{a2} follows from \eqref{a4} and $Spec(A_{h})=Spec(B_{h})-1$. The proof of lemma \ref{lemA1} is
complete.
\ep

\begin{lem}\label{lemcoho}
Let $\NN=\NN_{1}\oplus...\oplus \NN_{\r} $ be the free up to rank $\r$
nilpotent Lie algebra with $p$ generators. Let $(Y_{1},...,Y_{p})$ be a basis of
$\NN_{1}$ and let $(\mathcal Z_{1},...,\mathcal Z_{p})$ be the right invariant vector fields on $\NN$
such that $\mathcal Z_{j}(0)=Y_{j}$.
Let $\mathcal S(\NN)$ be the Schwartz space of $\NN$. Let $\varphi \in \mathcal S(\NN)$,  be such that $\int_{\NN} \varphi dx =0$.
Then there exists $\varphi_{k}\in \mathcal S(\NN)$ such that
\be\label{a5c}
\varphi= \sum_{k=1}^p \mathcal Z_{k}(\varphi_{k})
\ee
\end{lem}

\bp
Let $Y^\alpha= H_{\alpha}(Y_{1},...,Y_{p})$ and let $\mathcal Z^\alpha$ be the right invariant vector fields on $\NN$
such that $\mathcal Z^\alpha(0)=Y^\alpha$. Let $u_{\alpha}, \alpha\in \mathcal A$ be the coordinates on $\NN$
associated to the basis $(Y^\alpha, \alpha\in \mathcal A)$ of $\NN$. Let  $\partial_{\alpha}$
be the derivative in the direction of $u_{\alpha}$. Let $\varphi \in \mathcal S(\NN)$
such that $\int_{\NN}\varphi dx=0$.
Using the Fourier transform in coordinates $(u_{\alpha})$, and $\hat\varphi (0)=0$, one get easily that there
exists functions $\psi_{\alpha}\in \mathcal S(\NN)$ such that

\be\label{a5}
\varphi= \sum_{\alpha\in \mathcal A} \partial_{\alpha}(\psi_{\alpha})
\ee

By \eqref{3.10}, the vector field $\mathcal Z^\alpha$ is of the form
$$\mathcal Z^\alpha= \partial_{\alpha}+ \sum_{\vert \beta \vert>\vert\alpha\vert}
p_{\alpha,\beta}(u_{<\vert\beta\vert})\ \partial_{\beta}=\partial_{\alpha}+ \sum_{\vert \beta \vert>\vert\alpha\vert}
\partial_{\beta}\  p_{\alpha,\beta}(u_{<\vert\beta\vert})$$
where the $ p_{\alpha,\beta}$ are polynomials in $u$ depending only on
$(u_{1},...,u_{j})$ with $j<\vert\beta\vert$. Therefore, there exists polynomials 
$q_{\alpha,\beta}$ such that 

$$\partial_{\alpha}= \mathcal Z^\alpha+\sum_{\vert \beta \vert>\vert\alpha\vert}
\mathcal Z^\beta \  q_{\alpha,\beta} $$
Since the Schwartz space $\mathcal S(\NN)$ is stable by multiplication by polynomials, 
we get from \eqref{a5} that there exists $\phi_{\alpha}\in \mathcal S(\NN)$ such that

\be\label{a6}
\varphi= \sum_{\alpha\in \mathcal A} \mathcal Z^\alpha (\phi_{\alpha})
\ee
For $\vert\alpha\vert>1$, there exists $j\in \{1,...,p\}$ and $\beta$
with $\vert\beta\vert=\vert\alpha\vert-1$ such that 
$\mathcal Z^\alpha=\mathcal Z_{j}\mathcal Z^\beta-\mathcal Z^\beta
\mathcal Z_{j}$. By induction on $\vert\alpha\vert$, since the Schwartz space $\mathcal S(\NN)$ is stable by the vector fields 
$\mathcal Z_{j}$,  this shows that for any
$\alpha$ and $\phi\in\mathcal S(\NN)$, there exists $\phi_{j}\in\mathcal S(\NN)$
such that $\mathcal Z^\alpha(\phi)=\sum_{j=1}^p \mathcal Z_{j}(\phi_{j})$.
Thus \eqref{a5c} follows from \eqref{a6}.
The proof of lemma \ref{lemcoho} is
complete.

\ep

\bibliography{hypo.bib}

\end{document}